%% file: Factorization_of_Coefficients_for_exp_and_log_in_Function_Fields.tex
\documentclass{amsart}

\usepackage{amssymb}

\usepackage{graphics}
\usepackage{tikz, tikz-cd}
\usetikzlibrary{matrix, arrows} 

\usepackage{verbatim}
\usepackage{hyperref}
\usepackage{enumitem}
\usepackage{url} 
\usepackage{subfiles} 

\newtheorem{thm}{Theorem}[section]
\newtheorem{prop}[thm]{Proposition}
\newtheorem{lem}[thm]{Lemma}
\newtheorem{cor}[thm]{Corollary}
\newtheorem*{thm*}{Theorem}
\newtheorem*{lem*}{Lemma}
\newtheorem*{cor*}{Corollary}
\newtheorem{conj}{Conjecture}[section]

\newtheoremstyle{named}{}{}{\itshape}{}{\bfseries}{.}{.5em}{\thmnote{#3}}
\theoremstyle{named}
\newtheorem*{namedthm*}{Theorem}

\theoremstyle{definition}
\newtheorem{defn}[thm]{Definition}

\theoremstyle{remark}
\newtheorem{rmk}[thm]{Remark}
\newtheorem*{rmk*}{Remark}

\numberwithin{equation}{section}

\newcommand\lrp[1]{\left({#1}\right)} 
\newcommand\bb[1]{\mathbb{#1}} 
\newcommand\mf[1]{\mathfrak{#1}} 
\newcommand\mc[1]{\mathcal{#1}} 

\def\e{\varepsilon} 

\DeclareMathOperator{\Gal}{Gal}
\DeclareMathOperator{\Spec}{Spec}

\DeclareMathOperator{\ord}{ord}

\DeclareMathOperator{\sgn}{sgn}

\begin{document}


\title[Factorization of exp and log coefficients in Function Fields]{Factorization of Coefficients for Exponential and Logarithm in Function Fields}


\author[Kwun Chung]{Kwun ``Angus'' Chung}
\address{Department of Mathematics, University of Michigan, Ann Arbor, MI 48109}
\email{angusck@umich.edu}





\begin{abstract}
Let $X$ be an elliptic curve or a ramifying hyperelliptic curve over $\bb F_q$. We will discuss how to factorize the coefficients of the exponential and logarithm series for a Hayes module over such a curve. This allows us to obtain $v$-adic convergence results for such exponential and logarithm series, for $v$ a ``finite'' prime. As an application, we can show that the $v$-adic Goss $L$-value $L_v(1,\Psi)$ is log-algebraic for suitable characters $\Psi$. 
\end{abstract}

\date{October 24, 2020}
\keywords{Hayes modules, shtuka functions, $v$-adic zeta values, $v$-adic $L$-values}


\maketitle

\setcounter{tocdepth}{1}
\tableofcontents{}

\section{Introduction}
\label{sec:Introduction}
\subfile{Sec1.tex}


\section{Notations and Backgound}
\label{sec:DrinfeldModules}

\subfile{Sec2.tex}
\section{Logarithm series for our curves}
\label{sec:LogSeriesExpression}
\subfile{Sec3.tex}

\section{Factorization by Divisors}
\label{sec:FactnByDivisors}

\subfile{Sec4.tex}


\section{\texorpdfstring{$v$}{v}-adic valuation for the exponential and logarithm coefficients}
\label{sec:vForCoeff}

\subfile{Sec5.tex}


\section{Application: Log-algebracity of \texorpdfstring{$v$}{v}-adic Goss \texorpdfstring{$L$}{L}-value}
\label{sec:Applications}

\subfile{Sec6.tex}


\section{Examples}
\label{sec:Eg}

\subfile{Sec7.tex}


\section{Possible directions for generalizations}
\label{sec:PossibleGeneralizations}
\subfile{Sec8.tex}




\bibliographystyle{amsalpha}
\bibliography{bibliography.bib}

\end{document}

%% file: Sec1.tex
In the number field case, it is known that the exponential series 
\[
e^z = \sum_{n = 0}^\infty \frac{z^n}{n!}
\]
converges for all $z \in \bb C_p$ with $|z|_p < p^{-\frac{1}{p-1}}$, and the logarithm series 
\[
\log(1+z) = \sum_{n = 1}^\infty (-1)^{n+1} \frac{z^n}{n}
\]
converges for all $z \in \bb C_p$ with $|z|_p < 1$. See \cite[\S IV.2]{Koblitz1984_book} for example. 

In the function field case, characteristic $p$ exponential and logarithm functions arising from the Drinfeld modules were first studied by Carlitz \cite{Carlitz1935} and Drinfeld \cite{Drinfeld1974}. Let $X$ be a smooth, projective, geometrically connected curve over $\bb F_q$. We pick a closed point $\infty$ on $X$, and let $\bb A := \Gamma(X - \infty, \mc O_X)$, $\bb K := \mathrm{Frac}(\bb A)$. The idea of a Drinfeld module is to define another $\bb A$-action (other than scalar multiplication) on the function fields and their field extensions. To distinguish the action and the space being act on, we set $A,K$ to be isomorphic copies of $\bb A,\bb K$ respectively, and use the non-bold characters to represent spaces being act on. 

We are particularly interested in a special class of Drinfeld modules, called $\sgn$-normalized Drinfeld modules, or Hayes modules. They were first investigated by Hayes \cite{Hayes1974}, \cite{Hayes1979}, who utilize them to work out the class field theory of function fields. The simplest case of a Hayes module is when $X = \bb P^1$, $\bb A = \bb F_q[t]$, $A = \bb F_q[\theta]$. What we get then is the usual Carlitz module \cite{Carlitz1935}. In this case, 
\[
e_\rho(z) = \sum_{n=0}^\infty \frac{1}{(\theta^{q^n} - \theta) (\theta^{q^{n-1}}-\theta)^q \cdots (\theta^q-\theta)^{q^{n-1}}} z^{q^n},
\]
\[
\log_\rho(z) = \sum_{n = 0}^\infty \frac{1}{(\theta - \theta^q) \cdots (\theta - \theta^{q^n})} z^{q^n}.
\]
Similar to the number field case, there is a $v$-adic convergence result for these series, where $v$ is a place of $K = \bb F_q(\theta)$ corresponding to a prime ideal $\mf p$ in $A$. This can be found, for instance, in \cite{Anderson_Thakur1990} (where we take the special case of dimension 1 from the paper).
\begin{itemize}
\item $e_\rho(z)$ converges in $\bb C_v$ for all $z \in \bb C_v$ with $v(z) > \frac{1}{q^{\deg \mf p}-1}$;
\item $\log_\rho(z)$ converges in $\bb C_v$ for all $z \in \bb C_v$ with $v(z) > 0$.
\end{itemize}
We call the functions on $\bb C_v$ defined by these series $e_{v,\rho}(z)$ and $\log_{v,\rho}(z)$ respectively. 

Anderson \cite{Anderson1996} gave an application of the $v$-adic convergence, showing that the value of $v$-adic $L$-function $L_v(1,\chi)$ for a character $A \to A/\mf p \to \bb C_\infty$ is ``log-algebraic''. Namely, it is of the form 
\[
\sum_{i=1}^s \alpha_i \log_{v,\rho} S_i
\]
for $\alpha_i, S_i \in \overline{K}$. A similar application is given by Anderson and Thakur \cite{Anderson_Thakur1990}, which uses a higher-dimensional version of the $v$-adic convergence to compute certain $v$-adic zeta values. 


\subsection{Main results}

In this paper, we are going to generalize the $v$-adic convergence result to other Hayes modules. Specifically, we study the case of an elliptic curve or a {\em ramifying hyperelliptic curve}, i.e. a curve with an affine model
\[
\Gamma(X - \infty, \mc O_X) = \bb F_q[t,y] / (y^2 + F_2(t)y - F_1(t))
\]
for some rational point $\infty$ on $X$. Fix a Hayes module $\rho$ for the pair $(X,\infty)$ and a place $v$ other than $\infty$ on such a curve. The exponential and logarithm series $e_\rho(z)$ and $\log_\rho(z)$ for $\rho$ have coefficients in a finite Galois extension of $K$ (see discussion in section \ref{sec:DrinfeldModules} before Theorem \ref{thm:expByShtFncThakur}). Fix an embedding $\overline K \to \bb C_v$. Our main result is that:
\begin{namedthm*}[Theorem \ref{thm:vadicConvOfExp}]
The power series $e_\rho(z)$ converges in $\bb C_v$ for $z \in \bb C_v$ with 
\[
v(z) > C_1 + C_2 \frac{1}{q^{C_3}-1},
\]
where $C_1, C_2, C_3$ are explicit constants, to be given in section \ref{sec:vForCoeff}. In particular, $C_1$ is related to the Drinfeld divisor $V$ (to be defined in section \ref{sec:DrinfeldModules}), $C_2$ to some ramification indices, and $C_3$ to the degree of $v$ and some inertial degrees. 
\end{namedthm*}
Moreover,  
\begin{namedthm*}[Theorem \ref{thm:vadicConvOfLog}]
The power series $\log_\rho(z)$ converges in $\bb C_v$ for $z \in \bb C_v$ with $v(z) > 0$.
\end{namedthm*}

Similar to the Carlitz case, we call the functions on $\bb C_v$ defined by these series $e_{v,\rho}(z)$ and $\log_{v,\rho}(z)$ respectively. 

As an application, we will prove a similar log-algebraicity theorem for $v$-adic $L$-values over these curves. 
\begin{namedthm*}[Theorem \ref{thm:Lv(1,Psi)LogAlg}]
Let $\Psi$ be a multiplicative character of conductor $v$ on the group of $\bb A$-fractional ideals prime to $v$. Then $L_v(1,\Psi)$ is log-algebraic, i.e. there exists 
\[ \alpha_1, \cdots, \alpha_s, S_1, \cdots, S_s \in \overline K, \]
with $v(S_i) > 0$, such that 
\[
L_v(1,\Psi) = \sum_i \alpha_i \log_{v,\rho} S_i.
\]
\end{namedthm*}

To prove Theorems \ref{thm:vadicConvOfExp} and \ref{thm:vadicConvOfLog}, we recall an expression of $e_\rho(z)$ and $\log_\rho(z)$ given by Thakur, Anderson \cite[0.3.6, 0.3.8]{Thakur1993_ShtukasAJ}, Green and Papanikolas \cite[Corollary 3.5]{GreenPapanikolas2018}. These expressions are given in terms of a special function $f$, called the shtuka function, which was first studied by Thakur \cite{Thakur1993_ShtukasAJ} and realizes the correspondence that Drinfeld showed in \cite{Drinfeld1977}. The expressions are:
\[
e_\rho(z) = \sum_{n = 0}^\infty \left. \frac{1}{ff^{(1)} \cdots f^{(n-1)}} \right|_{\Xi^{(n)}} z^{q^n},
\]
\[
\log_\rho(z) = \sum_{n=0}^\infty \left. \frac{\delta^{(n+1)}}{\delta^{(1)} f^{(1)} \cdots f^{(n)}} \right|_{\Xi} z^{q^n},
\]
where $\Xi$ is the point corresponding to $\bb A \to A \to \bb C_\infty$, and $\delta$ is a function derived from the shtuka function $f$, and $(n)$ is the Frobenius twist (see section \ref{sec:DrinfeldModules}). 

We then proceed to study the coefficients by a factorization. The factorization will be done via Proposition \ref{prop:factorizationByDivisors}. Essentially, it tells us that the ideal generated by a sufficiently integral function can be factorized into a product of ideals given by the zeros of the function. Applying this to our coefficients, we will be able to obtain the factorization of the exp and log coefficients from the divisors of the shtuka function and its twists. We then use some basic theory of extension of valuations to study the valuations of each of the factors in section \ref{sec:vForCoeff}, thus proving theorems \ref{thm:vadicConvOfExp} and \ref{thm:vadicConvOfLog}.


\subsection{Outline of paper}

We will first recall some facts about Drinfeld modules, Hayes modules, and shtuka functions in section \ref{sec:DrinfeldModules}. In section \ref{sec:LogSeriesExpression}, we will derive a simplified expression for the logarithm series for our curves, following Green and Papanikolas \cite[Corollary 3.5]{GreenPapanikolas2018}. 
Section \ref{sec:FactnByDivisors} is where we prove the proposition \ref{prop:factorizationByDivisors} mentioned above that allows us to factorize functions according to its divisors. We will then prove theorems \ref{thm:vadicConvOfExp} and \ref{thm:vadicConvOfLog} in section \ref{sec:vForCoeff}. As an application, we will show the log-algebraicity theorem \ref{thm:Lv(1,Psi)LogAlg} for the $v$-adic Goss $L$-value $L_v(1,\Psi)$ in section \ref{sec:Applications}. Then we will give some examples in section \ref{sec:Eg}, as well as discussing possible generalizations in section \ref{sec:PossibleGeneralizations}.


\subsection{Acknowledgment}
The author would like to thank his advisor Kartik Prasanna for his continued guidance and support, and Dinesh Thakur for his helpful comments for this paper.

%% file: Sec2.tex
\subsection{Notation and Definition}

Let $X$ be a smooth, projective, geometrically connected curve of genus $g$ over $\bb F_q$. We assume there is a rational point $\infty \in X(\bb F_q)$ such that $\bb A := \Gamma(X - \infty, \mc O_X)$ has a model 
\[
\bb F_q[t,y] / (y^2 + F_2(t)y - F_1(t)),
\]
with $F_1,F_2 \in \bb F_q[t]$, $F_1$ monic of degree $2g+1$, $F_2$ of degree at most $g$. When $g = 1$, this is simply the usual Weierstrass model for elliptic curves. And when $g \geq 2$ this is what we call a {\em ramifying hyperelliptic curve}. The word ``ramifying'' signifies that the place $\infty$ ramifies in the extension $\mathrm{Frac}(\bb A) / \bb F_q(t)$. We will denote by $F(t,y)$ the polynomial $y^2+ F_2(t) y - F_1(t)$. 

We set $A$ to be a ring isomorphic to $\bb A$, but with $\theta, \eta$ replacing $t,y$ respectively. This is to distinguish the non-trivial action of $\bb A$ coming from the Drinfeld module structure from the scalar multiplication by $A$. For $a \in \bb A$, we set $a|_{\Xi}$ to be the corresponding element in $A$. Here is a table of notations. 

\begin{center}
\begin{tabular}{c | l}
$K$, $\bb K$ & $\mathrm{Frac} (A)$, $\mathrm{Frac} (\bb A)$ \\
$\bb H$ & Hilbert class field of $\bb K$, splitting completely at $\infty$ \\
$\bb K_\infty$ & completion of $\bb K$ at $\infty$ \\
$H, K_\infty$ & extensions of $K$, isomorphic to $\bb H, \bb K_\infty$ respectively \\
$\mc O_{K'}$ & ring of integers in a field $K'$ (can be local or global) \\
$\bb C_\infty$ & completion of $\overline{K_\infty}$ \\
$\bb C_v$ & similar to $\bb C_\infty$, but using a finite place $v$ instead. \\
$\overline X$ & $X_{\bb C_\infty}$ \\
$\Xi$ & closed point of $\overline X(\bb C_\infty)$ corresponding to $\bb A \to A \to \bb C_\infty$ \\
$\overline \infty$ & closed point of $\overline X(\bb C_\infty)$ above $\infty \in X(\bb F_q)$ 
\end{tabular}
\end{center}
In general, bold characters (except $\bb C_\infty$ and $\bb C_v$) represent spaces with variables, and non-bold characters (and $\bb C_\infty$, $\bb C_v$) represent spaces with scalars. 

For $P$ a point on $\overline X$, its {\em Frobenius twist} $P^{(1)}$ is the point obtained by taking $q$-th power on all coordinates of $P$ (in any affine model of $\overline X$ containing $P$). This induces the Frobenius twist map on the group of divisors, as well as on the function field of $\overline X$, which we also denote by $- \mapsto -^{(1)}$. For example, we have $(ct)^{(1)} = c^q t$ for $c \in \bb C_\infty$. For a differential $\omega = f dg$ on $\overline X$, we set $\omega^{(1)} := f^{(1)} dg^{(1)}$. Finally, the notation $- \mapsto -^{(n)}$ stands for applying the twist $n$ times.

In this setup, the residue degree of $\infty$ is $d_\infty = 1$. We define the degree function on $\bb A$ by $\deg a = -v_\infty(a)$. In particular, $\deg t = 2$ and $\deg y = 2g+1$. 

We need an analogous notion of ``leading coefficient'' for the ring $\bb A$. Since $\bb A$ has an $\bb F_q$-basis $\{t^i, yt^i\}_{i \geq 0}$, the notion of leading coefficient can be taken to be the coefficient of the term with highest degree. We define a sign function on $\bb A$ by setting $\sgn a$ to be the leading coefficient of $a$ with respect to the basis, and $\sgn 0 = 0$. In particular $\sgn t = \sgn y = 1$. For a detailed definition of a general sign function, see \cite[\S7.2]{Goss1998_BasicStrcBook} or \cite[1.1.1]{Thakur2004}. Since $d_\infty = 1$, other sign functions are just $\bb F_q^\times$ multiples of this. We will stick to this sign function for simplicity, but everything in this paper works for any sign function we pick. 

The ring $\bb C_\infty \otimes_{\bb F_q} \bb A$ has the same basis $\{t^i, yt^i\}_{i \geq 0}$ as a $\bb C_\infty$-module. We extend the sign function $\sgn$ to a function 
\[
\widetilde\sgn: \bb C_\infty \otimes_{\bb F_q} \bb A \to \bb C_\infty,
\]
by using the leading coefficient with respect to the same basis as in the definition of $\sgn$. 

\subsection{Hayes modules}

Our main objects of study are Hayes modules. They should be considered as ``good representatives'' of isomorphism classes of rank 1 Drinfeld modules. We will recall the definition and some background here. Recall that $\bb C_\infty\{\tau\}$ is the ring of twisted polynomials over $\bb C_\infty$ (see \cite[\S1]{Goss1998_BasicStrcBook}).
\begin{defn}
A {\em $\sgn$-normalized Drinfeld module}, or equivalently a {\em Hayes module}, is an $\bb F_q$-algebra homomorphism 
\[
\rho: \bb A \to \bb C_\infty \{\tau\}
\]
such that 
\[
\rho_a = a|_{\Xi} + \cdots + \sgn(a) \tau^{\deg a}
\]
\end{defn}

Such objects can be constructed, following the procedure provided by Thakur \cite{Thakur1993_ShtukasAJ}. By the work of Drinfeld \cite{Drinfeld1974}, there is an effective divisor $V$ on $\overline X$ such that the divisor
\[
V^{(1)} - V + (\Xi) - (\overline \infty)
\]
is principal. Let $f$ be the unique function in the function field of $\overline X$ such that 
\[
\mathrm{div} (f) = V^{(1)} - V + (\Xi) - (\overline \infty),
\]
and such that $\widetilde\sgn (f) = 1$. This is called a {\em shtuka function} for the triple $(V, \overline \infty, \sgn)$. 

Thakur \cite{Thakur1993_ShtukasAJ} observed that the set 
\[
1, f, ff^{(1)}, ff^{(1)} f^{(2)}, \cdots
\]
is a $\bb C_\infty$-basis for the vector space $\Gamma(\overline X - \infty, \mc O_{\overline X}(V))$. For 
\[
a \in \bb A \hookrightarrow \Gamma(\overline X-\infty, \mc O_{\overline X}(-V)),
\]
we can write 
\[
a = \rho_{a,0} + \rho_{a,1} f + \cdots + \rho_{a,\deg a} f f^{(1)} \cdots f^{(\deg a-1)}.
\]
where $\rho_{a,i} \in \bb C_\infty$. Define 
\[
\rho_a := \rho_{a,0} + \rho_{a,1} \tau + \cdots + \rho_{a,\deg a} \tau^{\deg a}.
\]
Then $\rho$ is a Hayes module for $\sgn$. 

\begin{rmk}
\label{rmk:normalizationOfShtukaFunction}
For a general curve with the residue degree at infinity $d_\infty$ bigger than 1, instead of $\widetilde\sgn (f) = 1$, which is the condition stated in \cite[0.3.5]{Thakur1993_ShtukasAJ} and \cite[\S7]{Goss1998_BasicStrcBook}, we should take 
\[
\widetilde\sgn (f f^{(1)} \cdots f^{(d_\infty-1)} ) = 1.
\]
Otherwise, $\rho$ will not have the correct top coefficient, i.e. it will not be $\sgn$-normalized. The main problem with $\widetilde\sgn (f) = 1$ is that we cannot determine what $\widetilde\sgn (f^{(1)})$ is. In fact $\widetilde\sgn (f^{(1)})$ can be transcendental over $\bb F_q$! See \S \ref{eg:sgnTranscendental} for an example with such a function. Generally, we can only calculate $\widetilde\sgn (F^{(d_\infty)})$ from $\widetilde\sgn (F)$, but not for any twists in between 0 and $d_\infty$. See \cite[(14)]{Thakur2004}. 

This condition 
\[
\widetilde\sgn (f f^{(1)} \cdots f^{(d_\infty-1)} ) = 1.
\]
does not determine a shtuka function uniquely. Instead, a shtuka function is defined up to a $\dfrac{q^{d_\infty}-1}{q-1}$-th root of unity from a triple $(V, \overline \infty, \sgn)$. This matches Goss's calculation \cite[7.2.18]{Goss1998_BasicStrcBook} on the number of Hayes modules for $\sgn$ within an isomorphism class. Thus, we can able to specify a unique shtuka function $f$ to a Hayes module $\sgn$, by requiring that 
\[
a = \rho_{a,0} + \rho_{a,1} f + \cdots + \rho_{a,\deg a} f f^{(1)} \cdots f^{(\deg a-1)}.
\]
\end{rmk}
~\\
Any Drinfeld module over $\bb C_\infty$ comes equipped with an exponential series and a logarithm series. They are power series in $z^q$ over $\bb C_\infty$, with coefficient of $z$ equals 1 and satisfies the following functional equations for all $a \in \bb A$:
\[
e_\rho(a|_{\Xi}z) = \rho_a (e_\rho(z)), \qquad \log_\rho(\rho_a(z)) = a|_{\Xi} \log_\rho(z).
\]
Since a Hayes module $\rho$ has coefficients in a finite Galois extension of $K$, namely the {\em field of definition} (see \cite[\S7.3]{Goss1998_BasicStrcBook}, \cite[\S3.3]{Thakur2004}), the coefficients of $e_\rho$ and $\log_\rho$ are in the same finite Galois extension, from the above functional equations. 

Thakur and Anderson \cite{Thakur1993_ShtukasAJ} gave a formula for the series in terms of the corresponding shtuka function. 
\begin{thm}
\label{thm:expByShtFncThakur}
\cite[0.3.6]{Thakur1993_ShtukasAJ} The exponential series for $\rho$ is given by 
\[
e_\rho(z) = z + \sum_{n = 1}^\infty \left. \frac{1}{f f^{(1)} \cdots f^{(n-1)}} \right|_{\Xi^{(n)}} z^{q^n}.
\]
\end{thm}

To describe the logarithm, we need a differential 
\[
\omega \in \Gamma(X, \Omega^1(-V + 2(\overline \infty)))
\]
such that 
\[
\mathrm{Res}_{\Xi} \frac{\omega^{(1)}}{f} = 1.
\]

\begin{thm}
\label{thm:logByResidue}
\cite[0.3.8]{Thakur1993_ShtukasAJ} The logarithm series for $\rho$ is given by 
\[
\log_\rho(z) = \sum_{n = 0}^\infty \mathrm{Res}_{\Xi} \frac{\omega^{(n+1)}}{f f^{(1)} \cdots f^{(n)}} z^{q^n}.
\]
\end{thm}

\begin{rmk}
\label{rmk:NormalizationOfDifferential}
For a general curve with $d_\infty > 1$, we need to be more careful about the poles of $\omega$. Now, $\overline \infty$ is a fixed closed point in $\overline X(\bb C_\infty)$ above $\infty$. Suppose we pick our shtuka function such that 
\[
\mathrm{div} (f) = V^{(1)} - V + (\Xi) - (\overline \infty^{(1)}).
\]
The definition stated in \cite[0.3.7]{Thakur1993_ShtukasAJ} and \cite[\S7.11]{Goss1998_BasicStrcBook} specifies that we take $\Omega^1(V - 2(\overline \infty))$. When $d_\infty > 1$, the differential
\[
\frac{\omega^{(k+1)}}{f f^{(1)} \cdots f^{(k)}}
\]
will have a simple pole at $\overline \infty^{(k+1)}$ for $1 \leq k < d_\infty$, which is undesired. Thakur has fixed this by picking the differential $\omega$ from $\Omega^1(V - (\overline \infty) - (\overline \infty^{(-1)}))$ instead. See \cite[(14)]{Thakur_update}. 

We will use the corrected normalizations in the previous remark \ref{rmk:normalizationOfShtukaFunction} and this one to calculate an example in the case $d_\infty > 1$ in section \ref{sec:Eg}. 
\end{rmk}

%% file: Sec3.tex
For elliptic curves, Green and Papanikolas \cite[Corollary 3.5]{GreenPapanikolas2018} give a simplified expression for the logarithm series. We will show that such an expression also exists for ramifying hyperelliptic curves. To formulate this expression, we first need to write down the shtuka function $f$ and the differential $\omega$ more explicitly. 

\begin{prop}
There exist polynomials $\delta, Q \in \mc O_H[t]$ of degree $g$ in $t$ such that: 

\begin{enumerate}
\item $\delta$ is monic;

\item setting
\[
\nu = y - \eta - Q(t),
\]
the shtuka function $f$ has a presentation 
\[
f = \frac{\nu}{\delta};
\]

\item $\widetilde\sgn (\nu) = \widetilde\sgn (\delta) = 1$;

\item in the affine model $\overline X - \overline \infty$ with coordinate ring $\bb C_\infty \otimes_{\bb F_q} \bb A$, the coordinates of the zeros of $\nu$ and $\delta$ are integral over $\mc O_H$.

\end{enumerate}
\end{prop}

The proof is by long division, as in the proof of \cite[Theorem V]{Thakur1992_DMArith}. 

\begin{proof}
The shtuka correspondence from section \ref{sec:DrinfeldModules} tells us that 
\begin{align*}
t &= \theta + \rho_{t,1} f + f f^{(1)}, \\
y &= \eta + \rho_{y,1} f + \rho_{y,2} f f^{(1)} + \cdots + f f^{(1)} \cdots f^{(2g)}.
\end{align*}
Rearranging, we get 
\begin{align}
(\theta - t) + \rho_{t,1} f + f f^{(1)} = 0, \tag{1} \label{eq:LongDiv1} \\
(\eta - y) + \rho_{y,1} f + \rho_{y,2} f f^{(1)} + \cdots + f f^{(1)} \cdots f^{(2g)} = 0. \tag{2} \label{eq:LongDiv2}
\end{align}

The long division starts as follows: twist equation (\ref{eq:LongDiv1}) by $2g-1$ times, multiply by $-f f^{(1)} \cdots f^{(2g-2)}$, then add (\ref{eq:LongDiv2}), i.e. we have $(\ref{eq:LongDiv2}) - f f^{(1)} \cdots f^{(2g-2)} \cdot ((\ref{eq:LongDiv1})^{(2g-1)})$: 
\begin{multline}
(\eta - y) + \cdots + \rho_{y,2g-2} f f^{(1)} \cdots f^{(2g-3)} \\
 + (\rho_{y,2g-1} - (\theta - t)) f f^{(1)} \cdots f^{(2g-2)} + (\rho_{y,2g} - \rho_{t,1}) f f^{(1)} \cdots f^{(2g-1)} = 0.
\tag{3} \label{eq:LongDiv3}
\end{multline}
Continue with process of long division. 

{\bf Claim: } For $3 \leq n \leq 2g+1$, after $n-2$ steps in the long division, equation ($n$) is of the form
\begin{multline}
(\eta - y) + \cdots + \rho_{y,2g-n+1} f f^{(1)} \cdots f^{(2g-n)} \\
+ P_n(t) f f^{(1)} \cdots f^{(2g-n+1)} + Q_n(t) f f^{(1)} \cdots f^{(2g-n+2)} = 0, \tag{$n$}
\end{multline}
where $P_n(t), Q_n(t)$ are polynomials in $\mc O_H[t]$ of degree (in $t$)
\[
\deg_t P_n(t) = \left\lfloor \frac{n-1}{2} \right\rfloor, \quad
\deg_t Q_n(t) = \left\lfloor \frac{n-2}{2} \right\rfloor. 
\]
Moreover, $P_n$ is monic when $n$ is odd, and $Q_n$ is monic when $n$ is even. \\

This claim can be seen by induction. In particular, the inductive step is by dividing equation ($n-1$) by equation (\ref{eq:LongDiv1}), which is very similar to how we obtain equation (\ref{eq:LongDiv3}) from equations (\ref{eq:LongDiv1}) and (\ref{eq:LongDiv2}). 

Now, take equation ($2g+1$) and divide it by equation (\ref{eq:LongDiv1}) to obtain 
\[
(\eta - y - (\theta - t) Q_{2g+1}(t)) + (P_{2g+1}(t) - \rho_{t,1} Q_{2g+1}(t)) f = 0.
\]
Define 
\[
\delta := P_{2g+1}(t) - \rho_{t,1} Q_{2g+1}(t), \quad
Q := (t-\theta) Q_{2g+1}(t). 
\]

From the claim, $\deg_t P_{2g+1}(t) = g$, $\deg_t Q_{2g+1}(t) = g-1$, and $P_{2g+1}$ is monic. Thus,
\begin{enumerate}
\item $\delta$ and $Q$ are in $\mc O_H[t]$ and have degree $g$ in $t$;
\item $\delta$ is monic;
\item by setting $\nu := y - \eta - Q(t)$, we have $\displaystyle f = \frac{\nu}{\delta}$;
\item $\widetilde\sgn (\delta) = 1$ comes directly from $\delta$ being a monic polynomial in $t$;
\item since the degrees of $y$ and $Q(t)$ in $\bb A$ are $2g+1$ and $2g$ respectively, $\widetilde\sgn (\nu) = 1$. 
\end{enumerate}

What remains to be checked is the integrality of the zeros in the affine model. Since 
\[
\mathrm{div} (f) = V^{(1)} - V + (\Xi) - (\overline \infty),
\]
we must have $\mathrm{div} (\delta) \geq V$. Set $V' = \mathrm{div} \delta - V + 2g(\overline \infty)$, that is, 
\[
\mathrm{div} (\delta) = V + V' - 2g(\overline \infty).
\]
Since $\delta$ is a monic polynomial in $t$ with degree $g$ in $t$, $V'$ is the effective divisor of degree $g$. Then 
\[
\mathrm{div} (\nu) = V^{(1)} + V' + (\Xi) - (2g+1)(\overline \infty).
\]
Thus the zeros of $\delta$ and $\nu$ are given by $V^{(1)}$, $V$, $V'$, and $\Xi$. By definition $\Xi$ has coordinates $(\theta, \eta)$, both are in $\mc O_H$. Also the coordinates for $V^{(1)}$ are just $q$-th power of the coordinates for $V$. It remains to show that all points in the support of $V$ and $V'$ have coordinates integral over $\mc O_H$. 

The $t$-coordinates of all points in the support of $V$ and $V'$ are precisely the the roots of the polynomial $\delta$, when viewed as a polynomial in $t$. Since $\delta$ is monic and has coefficients in $\mc O_H$, the $t$-coordinates are integral. 

Plug in the $t$-coordinates for $V$ and $V'$, all of which integral over $\mc O_H$, into $F(t,y) = y^2 + F_2(t) y - F_1(t)$. We can see each of the $y$-coordinates for $V$ and $V'$ as a root of a monic quadratic polynomial over some integral extension of $\mc O_H$. Therefore the $y$-coordinates for $V$ and $V'$ are also integral. 

\end{proof}

With this presentation of the shtuka function, we can write down the differential $\omega \in H^0(\overline X, \Omega_{\overline X}^1(- V + 2(\overline \infty)))$ as in the previous subsection explicitly. This expression of the differential and the expression of the logarithm series derived from it (proposition \ref{prop:logByShtuka}) below is generalizes a result of Green and Papanikolas \cite[Corollary 3.5]{GreenPapanikolas2018} to our setting. 

\begin{prop}
Set 
\[
\omega := \frac{\delta}{2y + F_2(t)} dt. 
\]
Then $\omega \in H^0(\overline X, \Omega_{\overline X}^1(- V + 2(\overline \infty)))$, and
\[
\mathrm{Res}_{\overline \infty} \frac{\omega^{(1)}}{f} = -1. 
\]
\end{prop}

\begin{proof}
We first analyze the divisor of the differential $\dfrac{dt}{2y + F_2(t)}$. The function $t$ is holomorphic everywhere except at $\overline \infty$. Thus $dt$ cannot have pole except at $\overline \infty$. Recall that $v_\infty(t) = -2$ and $v_\infty(y) = -(2g+1)$, so we can pick $u := \dfrac{t^g}{y}$ to be a uniformizer at $\overline \infty$. First, we differentiate $F(t,y) = 0$ and see that 
\[
(2y + F_2(t))dy = (F_1'(t) - F_2'(t)y) dt.
\]
Next, we analyze $du$:
\begin{align*}
du &= \frac{gt^{g-1}}{y} dt - \frac{t^g}{y^2} dy \\
&= \frac{gt^{g-1}}{y} dt - \frac{t^g}{y^2} \frac{F_1'(t) - F_2'(t)y}{2y + F_2(t)} dt \\
&= \frac{t^{g-1}}{y^2(2y + F_2(t))} \left( g (2y^2 + F_2(t)y) - t(F_1'(t) - F_2'(t)y) \right) dt \\
&= \frac{t^{g-1}}{y^2(2y + F_2(t))} \left( g (2F_1(t) - F_2(t)y) - t(F_1'(t) - F_2'(t)y) \right) dt \\
&= \frac{t^{g-1}}{y^2(2y + F_2(t))} \left( (2gF_1(t) - tF_1'(t)) + (-gF_2(t)y + F_2'(t) ty) \right) dt. 
\end{align*}
Since $F_1(t)$ is monic of degree $2g+1$ in $t$, the polynomial $2gF_1(t) - tF_1'(t)$ is also of degree $2g+1$ in $t$, with leading coefficient $-1$. The term $-gF_2(t)y + F_2'(t) ty$ has fewer poles at $\overline \infty$ than $t^{2g+1}$ does. This shows that 
\[
\frac{dt}{2y + F_2(t)} = (-u^{2g-2} + O(u^{2g-1})) du.
\] 

Since $X$ is smooth, $dt$ and $dy$ cannot share any zero. By revisiting the equation 
\[
(2y + F_2(t))dy = (F_1'(t) - F_2'(t)y) dt,
\]
all zeros of $dt$ must also be zeros of $2y + F_2(t)$. That is, the differential $\dfrac{dt}{2y + F_2(t)}$ has no zeros except potentially at the poles of $2y+F_2(t)$, i.e. at $\overline \infty$. Since all differentials on $X$ have degree $2g-2$ and $\dfrac{dt}{2y+F_2(t)}$ has degree $2g-2$ at $\overline \infty$, the differential $\dfrac{dt}{2y + F_2(t)}$ has no zero nor pole away from $\overline \infty$. Therefore, 
\[
\omega = \frac{\delta}{2y + F_2(t)} dt
\]
has zeros at $V$, $V'$, a double pole at $\overline \infty$, and no other zeros or poles. In particular, $\omega \in H^0(\overline X, \Omega_{\overline X}^1(- V + 2(\overline \infty)))$. 

To prove the second statement about the residue, let us expand the differential 
\[
\frac{\omega^{(1)}}{f} = \frac{1}{f} \frac{\delta^{(1)}}{2y + F_2(t)} dt
\]
with respect to the uniformizer $u = \dfrac{t^g}{y}$. By definition of $f$ we have $\widetilde\sgn (f) = 1$. Since $\sgn t = \sgn y = 1$, this shows that 
\[
f = u^{-1} + O(u).
\]
Since $\delta^{(1)}$ is a monic polynomial of degree $g$ in $t$, we have that 
\[
\delta^{(1)} = u^{-2g} + O(u^{-2g+1}).
\]
Therefore, 
\[
\frac{\omega^{(1)}}{f} = -u^{-1} + O(1),
\]
showing that 
\[
\mathrm{Res}_{\overline \infty} \frac{\omega^{(1)}}{f} = -1.
\]
\end{proof}

With the explicit description of $\omega$, we can then write down formulae for the residues appearing in Theorem \ref{thm:logByResidue}. Note that $t - \theta$ is a uniformizer at $\Xi$, with $d(t-\theta) = dt$. Thus, 
\begin{align*}
\mathrm{Res}_{\Xi} \frac{\omega^{(n+1)}}{f f^{(1)} \cdots f^{(n)}}
&= \mathrm{Res}_{\Xi} \frac{\delta^{(n+1)}}{f f^{(1)} \cdots f^{(n)}} \frac{dt}{2y + F_2(t)} \\
&= \left. \frac{\delta^{(n+1)}}{\delta^{(1)} f^{(1)} \cdots f^{(n)}} \right|_{\Xi} \mathrm{Res}_{\Xi} \frac{\delta^{(1)} dt}{(2y + F_2(t))f} \\
&= \left. \frac{\delta^{(n+1)}}{\delta^{(1)} f^{(1)} \cdots f^{(n)}} \right|_{\Xi} \mathrm{Res}_{\Xi} \frac{\omega^{(1)}}{f} \\
&= \left. \frac{\delta^{(n+1)}}{\delta^{(1)} f^{(1)} \cdots f^{(n)}} \right|_{\Xi}.
\end{align*}

\begin{prop}
(cf. \cite[Cor 3.5]{GreenPapanikolas2018})
\label{prop:logByShtuka}
\[
\log_\rho(z) = \sum_{n =0}^\infty \left. \frac{\delta^{(n+1)}}{\delta^{(1)} f^{(1)} \cdots f^{(n)}} \right|_{\Xi} z^{q^n}. 
\]
\qed
\end{prop}

%% file: Sec4.tex
In this section, we will state and prove a proposition that factorizes a ``nice enough'' function $G$ according to divisors. We will do this more generally, so one can hope to apply this proposition to a general curve. 

Let $X / \bb F_q$ be a smooth projective geometrically connected curve over $\bb F_q$, $\infty$ a closed point (not necessarily rational), $\bb A := \Gamma(X-\infty, \mc O_X)$, and $\bb K$, $A$, etc. as before. Denote by $\bb F_\infty$ the residue field at $\infty$. Let $R$ be a finitely-generated integral domain over $\bb F_q$, 
$F = \mathrm{Frac} (R)$ and $L = \overline{F}$. All tensor products will be over $\bb F_q$ unless otherwise specified.

Since $\bb F_q$ is perfect, $X$ is geometrically normal. Thus $F \otimes \bb A$ is a domain, integrally closed in its field of fractions. In particular, $R \otimes \bb A$ is an integral domain. It is easy to check that
\[
\mathrm{Frac} (R \otimes \bb A) = \mathrm{Frac} (F \otimes \bb A),
\]
the ``$=$'' sign means canonical isomorphism. Note that the latter field is the function field of $X_F := X \times_{\bb F_q} \Spec F$. 

Our goal is to analyze $R \otimes \bb A$-ideals, using the information we can obtain when we view elements in $R \otimes \bb A$ as functions on $X_F$ or $X_L := X \times_{\bb F_q} \Spec L$. 


\subsection{Functions in the subring \texorpdfstring{$R \otimes \bb A$}{R tensor bold A}}

Suppose we have an element $G \in R \otimes \bb A$. We can view $G$ as a meromorphic function on the curve $X_F$ or $X_L$. In the latter case, since $L$ is algebraically closed, we can write 
\[
\mathrm{div} (G) = \sum_i Z_i - \sum_j P_j,
\]
where $Z_i, P_j \in X(L)$ are the zeros and poles of $G$ respectively. Since $G \in F \otimes \bb A = \Gamma(X_F - \infty, \mc O_X)$, it is regular away from $\infty$. That is, all poles of $G$ must lie above $\infty$. 


We now fix a model for $\bb A$. Suppose 
\[
\bb A = \bb F_q[t_1, t_2, \cdots, t_n] / (F_1, F_2, \cdots, F_m).
\]

With the model, we can talk about the zeros of $G$ in terms of coordinates. Say $Z_i$ is given by $t_1 = a_{i,1}, t_2 = a_{i,2}, \dots, t_n = a_{i,n}$, where $a_{i,k} \in L$. 

\begin{defn}
Suppose $G \in R \otimes \bb A$. If there exists a model for $\bb A$ such that for all $i,k$, we have $a_{i,k} \in R$ with $a_{i,k}$ the coordinates of the zeros of $G$, then we say that {\bf all zeros of $G$ are in $R$}. 
\end{defn}

If all zeros of $G$ are in $R$, we can reduce the zeros modulo primes from $R$. This is what we will use for the proof of the proposition of the next subsection. As a remark, it is easy to see that if all zeros of $G$ are in $R$ for one model of $\bb A$, then it is true for all models of $\bb A$. 


\subsection{Infinities and \texorpdfstring{$\sgn$}{sgn}}

Consider the points lying above $\infty \in X(\bb F_q)$ in the tower:
\begin{center}
\begin{tikzcd}
X_{L} \ar[dr, "i_1"] & & X_{\overline{\bb F_q}} \ar[dl, "i_2" above] \\
& X_{\bb F_\infty} \ar[d] & \\
& X = X_{\bb F_q} &
\end{tikzcd}
\end{center}

The point $\infty$, as an $\overline{\bb F_q}$ point in $X$, can only split or be inert in this tower, and it splits completely as long as the field of constants contains $\bb F_\infty$. Hence, the fibers of $\infty$ in $X_L$ and in $X_{\overline{\bb F_q}}$ are in natural bijection, given by 
\begin{align*}
\{\infty\} \times_{X} X_L &\stackrel{=}{\longrightarrow} \{\infty\} \times_X X_{\overline{\bb F_q}} \\
\overline{\infty} &\mapsto i_2^* (i_{1,*} (\overline \infty)).
\end{align*}
By abuse of notation, we use ``$\overline \infty$'' to denote both a point in $X_L$ above $\infty$, and the corresponding point in $X_{\overline{\bb F_q}}$. 

For a sign function $\sgn$ on $\bb K_\infty$, it can be extended to a function 
$$ \widetilde\sgn_{\overline \infty}: \mathrm{Frac}(L \otimes \bb A) \to L. $$ 
This extension depends on a choice of $\overline \infty \in X(L)$ above $\infty$, and is unique upon such a choice. See \cite[\S7.11]{Goss1998_BasicStrcBook}. In the same manner, we can also extend $\sgn$ to a function 
$$ \widetilde\sgn_{\overline \infty}: \mathrm{Frac}(\overline{\bb F_q} \otimes \bb A) \to \overline{\bb F_q}. $$
By abuse of notation, we denote by $\widetilde\sgn_{\overline \infty}$ both extensions with respect to $\overline \infty$. 


\subsection{Relating the ideal by zeros of the function}

We are now ready to state our proposition that factorizes ``sufficiently integral'' functions in terms of their divisor. Let $R[\bb F_\infty]$ to be the smallest subring of $L$ containing $R$ and $\bb F_\infty$. 

\begin{prop}
\label{prop:factorizationByDivisors}
Suppose we have $G \in R \otimes \bb A$ such that $\widetilde{\sgn}_{\overline \infty} (G) \in (R[\bb F_{\infty}])^\times$ for all choices of $\overline \infty$, and that all its zeros are in $R$. Fix a model for $\bb A$ and let $\{(a_{i,k})_k\}_i$ be the zeros of $G$. Then we have an equality of $R \otimes \bb A$-ideals
\[
(G) = \prod_i (t_1 - a_{i,1}, \cdots, t_n - a_{i,n}).
\]
\end{prop}

This is proved by reducing mod $\mf m$ for maximal ideals $\mf m$ in $R$. Let $\widetilde G := G \mod \mf m$. We first prove a lemma that computes the degree of $\widetilde G$. 
\begin{lem}
\label{lem:degRAfterRedn}
Fix $\mf m \subset R$ a maximal ideal, and let $k' := R/\mf m$, which is a finite extension of $\bb F_q$ since $R$ is finitely generated over $\bb F_q$. Fix a closed point $\overline \infty$ in $X_L$ above $\infty$, and we also denote by $\overline \infty$ the corresponding point on $X_{\overline{\bb F_q}}$. Then
\[
\ord_{\overline \infty} (G) \text{ in } X_L = \ord_{\overline \infty} (\widetilde G) \text{ in } X_{\overline{\bb F_q}}.
\]
\end{lem}

\begin{proof}
Fix a uniformizer $\pi$ in $\bb K$ of $\infty$. Then 
\[
G = f_{-m} \pi^{-m} + O(\pi^{-m+1}),
\]
where $f_{-m} \in F$. By assumption, $\widetilde{\sgn} G = f_{-m} (\sgn \pi)^{-m}$ is a unit in $R[\bb F_\infty]$, so $f_{-m}$ is a unit in $R[\bb F_\infty]$. In particular, $f_{-m}$ is in $F \cap R[\bb F_\infty] = R$, and hence $f_{-m} \in R^\times$ by going-up theorem. By reducing modulo $\mf m$, viewing the coefficients of the expansion of $G$ as in $F_{\mf m}$, we get 
\[
\widetilde{G} = \widetilde{f_{-m}} \pi^{-m} + O(\pi^{-m+1}),
\]
showing that $G$ and $\widetilde G$ has the same order of poles at $\overline \infty$. 

\end{proof}

\begin{proof}[Proof of Proposition \ref{prop:factorizationByDivisors}]
Let $M_1, M_2$ be the $R \otimes \bb A$-modules on the left and right side respectively. We will show that for all maximal ideals $\mf m$ of $R$, the $(R / \mf m) \otimes \bb A$-modules $\widetilde M_i := M_i / (\mf m \otimes \bb A) M_i$ are equal for $i = 1,2$. Then we will use Nakayama's lemma to conclude that $M_1 = M_2$. 



Since $G$ and $\widetilde G$ have the same number of poles at each (respective) $\overline \infty$, they have the same number of zeros. Thus in $X_{\overline{\bb F_q}}$, we can write 
\[
\mathrm{div} (\widetilde G) = \sum_i (t_1 - \widetilde{a_{i,1}}, \cdots, t_n - \widetilde{a_{i,n}}) - \sum \text{poles above } \infty.
\]

Since $X$ is a smooth curve over $\bb F_q$, the ring $\overline{\bb F_q} \otimes \bb A$ is noetherian (from $X$ being locally of finite type over $\bb F_q$), integrally closed (since $X$ is smooth, hence normal, hence geometrically normal as $\bb F_q$ is perfect), and has Krull dimension 1. That is, it is Dedekind. As a result, all non-zero ideals of $\overline{\bb F_q} \otimes \bb A$ can be factored into a product of maximal ideals. By the Nullstellensatz, maximal ideals of 
\[
\overline{\bb F_q}[t_1,\cdots,t_n]
\]
are of the form
\[
(t_1 - \alpha_1, \cdots, t_n - \alpha_n).
\]
Since ideals of $\overline{\bb F_q} \otimes \bb A$ are in natural bijection with ideals of $\overline{\bb F_q}[t_1,\cdots,t_n]$ which contain $F_1,\cdots,F_m$, we have as ideals in $\overline{\bb F_q} \otimes \bb A$, 
\[
(\widetilde G) = \prod_i (t_1 - \widetilde{a_{i,1}}, \cdots, t_n - \widetilde{a_{i,n}}).
\]
Intersecting the ideals with $k' \otimes \bb A$ (or equivalently, taking $\Gal(\overline{k'}/k')$-invariant), we have that $\widetilde{M_1} = \widetilde{M_2}$. 

We now proceed to show that $M_1 = M_2$. Let $M_3 = M_1 \cap M_2$ in $R \otimes \bb A$. From above, we know that 
\[
(\mf m \otimes \bb A)(M_i / M_3) = M_i / M_3
\]
for $i = 1,2$ and for all $\mf m \subset R$ maximal. Now for any maximal ideal $\mf M$ of $R \otimes \bb A$, the pullback of $\mf M$ to $R$ along $R \to R \otimes \bb A$ is maximal, because $R$ is of finite type over $\bb F_q$. Thus $\mf M$ contains $\mf m \otimes 1$ for some $\mf m \subset R$ maximal. As a result, 
\[
\mf M (M_i / M_3) = M_i / M_3
\]
for all maximal ideals $\mf M \subset R \otimes \bb A$. Since $M_i / M_3$ are finitely generated as $R \otimes \bb A$-modules, by Nakayama's lemma $M_i / M_3 = 0$. Therefore $M_1 = M_3 = M_2$. 

\end{proof}

The main way we apply this proposition is to evaluate at a certain closed point in $\overline X$ not above $\infty$. 

\begin{cor}
\label{cor:factorizationModXi}
Let $\xi$ be a closed point of $\overline X$ not above $\infty$, with coordinates in the model $\bb A$ given by $(\theta_1, \cdots, \theta_n)$ such that $\theta_k \in R$ for all $k$. Then as $R$-ideals, 
\[
(G|_{\xi}) = \prod_i (\theta_1 - a_{i,1}, \cdots, \theta_n - a_{i,n}).
\]
\qed
\end{cor}


\subsection{Applying the proposition to exp and log coefficients}

We analyze the log coefficients for elliptic curves and ramifying hyperelliptic curves in this subsection. Recall from section \ref{sec:DrinfeldModules} that the shtuka function has a presentation 
\[
f = \frac{\nu}{\delta},
\]
where $\nu, \delta \in \mc O_H$, $\widetilde\sgn (\nu) = \widetilde\sgn (\delta) = 1$, and all coordinates of their zeros are also integral over $\mc O_H$. Recall also that 
\[
\mathrm{div} (\nu) = V^{(1)} + V' + (\Xi) - (2g+1)(\overline \infty), \qquad \mathrm{div} (\delta) = V + V' - 2g(\overline \infty),
\]
where $V'$ is an effective divisor of degree $g$. Let $\{(t = \alpha_i, y = \beta_i)\}_{i=1}^g$ be the coordinates for $V$, and $\{(t = \alpha_i', y = \beta_i')\}_{i=1}^g$ be the coordinates for $V'$. Define $K_V$ to be the smallest field extension of $H$ containing all these coordinates, and $R$ the integral closure of $A$ in $K_V$. By applying propsition \ref{prop:factorizationByDivisors} to $\nu^{(k)}$, $\delta^{(k)}$ and the $R$ we just defined, we can see that as ideals of $R \otimes \bb A$,
\begin{align*}
(\nu^{(k)}) &= (t - \theta^{q^k}, y - \eta^{q^k}) \prod_i \left[ (t - \alpha_i^{q^{k+1}}, y - \beta_i^{q^{k+1}}) (t - \alpha_i'^{q^k}, y - \beta_i'^{q^k}) \right], \\
(\delta^{(k)}) &= \prod_i \left[ (t - \alpha_i^{q^k}, y - \beta_i^{q^k}) (t - \alpha_i'^{q^k}, y - \beta_i'^{q^k}) \right].
\end{align*}

Apply corollary \ref{cor:factorizationModXi} to the same set of data, with $\xi = \Xi$, we obtain that as $R$-ideals, 
\begin{align*}
(\nu^{(k)}|_{\Xi}) &= (\theta - \theta^{q^k}, \eta - \eta^{q^k}) \prod_i \left[ (\theta - \alpha_i^{q^{k+1}}, \eta - \beta_i^{q^{k+1}}) (\theta - \alpha_i'^{q^k}, \eta - \beta_i'^{q^k}) \right], \\
(\delta^{(k)}|_{\Xi}) &= \prod_i \left[ (\theta - \alpha_i^{q^k}, \eta - \beta_i^{q^k}) (\theta - \alpha_i'^{q^k}, \eta - \beta_i'^{q^k}) \right].
\end{align*}

Therefore, we have the following factorization of coefficients for $e_\rho$ and $\log_\rho$:

\begin{prop}
As fractional $R$-ideals, 
\begin{align*}
\lrp{ \left. \frac{1}{ff^{(1)} \cdots f^{(n-1)}} \right|_{\Xi^{(n)}} } 
&= \lrp{ \left. \frac{\delta \cdots \delta^{(n-1)}}{\nu \cdots \nu^{(n-1)}} \right|_{\Xi^{(n)}} } \\
&= \lrp{ \prod_{k=0}^{n-1} (\theta^{q^n} - \theta^{q^k}, \eta^{q^n} - \eta^{q^k})^{-1} } \\
&\phantom{==} \cdot \lrp{ \prod_i (\theta^{q^n} - \alpha_i, \eta^{q^n} - \beta_i) (\theta^{q^n} - \alpha_i^{q^n}, \eta^{q^n} - \beta_i^{q^n})^{-1} },
\end{align*}
and 
\begin{align*}
\left( \left. \frac{\delta^{(n+1)}}{\delta^{(1)} f^{(1)} \cdots f^{(n)}} \right|_{\Xi} \right) 
&= \left( \left. \frac{\delta^{(2)} \cdots \delta^{(n+1)}}{\nu^{(1)} \cdots \nu^{(n)}} \right|_{\Xi} \right) \\
&= \lrp{ \prod_{k=1}^n (\theta - \theta^{q^k}, \eta - \eta^{q^k})^{-1} } \\
&\phantom{==} \cdot \lrp{ \prod_i (\theta - \alpha_i'^{q^{n+1}}, \eta - \beta_i'^{q^{n+1}}) (\theta - \alpha_i'^q, \eta - \beta_i'^q)^{-1} }.
\end{align*}
\end{prop}


%% file: Sec5.tex

We will study the $v$-adic valuation for the coefficients of exponential and logarithm series for elliptic curves and ramifying hyperelliptic curves in this section, by looking at the $v$-adic valuation of each term in the factorization of ideals. For ease of notation, for $n \geq 0$, let $e_n,l_n$ be defined by 
\[
e_\rho(z) = \sum_{n \geq 0} e_n z^{q^n}, \qquad
\log_\rho(z) = \sum_{n \geq 0} l_n z^{q^n}.
\]


As a remark, none of the $l_n$ is 0 by \cite[Theorem 8.3.13]{Thakur2004}. 


\subsection{Notations from Elementary Number Theory}

Let us first fix some notation for this subsection. Recall that $\mf p \subset A$ is the prime corresponding to $v$. Suppose $\mf p$ has degree $d_{\mf p}$, i.e. $d_{\mf p} = -d_\infty v_\infty(\mf p) = -v_\infty(\mf p)$. Let $\mf p_\theta \in \bb F_q[\theta]$, $\mf p_\eta \in \bb F_q[\eta]$ be the monic generators of the ideals $\mf p \cap \bb F_q[\theta]$ and $\mf p \cap \bb F_q[\eta]$ respectively. By abuse of notation, we also use $\mf p_\theta$, $\mf p_\eta$ to denote the maximal ideal they generate in $\bb F_q[\theta]$ and $\bb F_q[\eta]$ respectively. By elementary number theory, we have that 
\[
\deg \mf p = \deg_\theta \mf p_{\theta} \cdot f_\theta = \deg_\eta \mf p_{\eta} \cdot f_\eta, 
\]
where $f_\theta$ is the inertial degree of $\mf p$ over $\mf p_{\theta}$, and similar for $f_\eta$. Note that $\deg \mf p$ is divisible by both $f_\theta$ and $f_\eta$, hence also by $\gcd(f_\theta, f_\eta)$. We also set $v_{\mf p_\theta}$ and $v_{\mf p_\eta}$ to be the valuation on $\bb F_q(\theta)$ and $\bb F_q(\eta)$ corresponding to $\mf p_\theta$ and $\mf p_\eta$ respectively, such that $v_{\mf p_\theta}(\mf p_\theta) = 1$ and $v_{\mf p_\eta}(\mf p_\eta) = 1$. Because we normalize $v$ so that the value group is $\bb Z$, for $b_1 \in \bb F_q(\theta)$ and $b_2 \in \bb F_q(\eta)$, viewing $b_1, b_2$ as functions in $K$, we have
\[
v(b_1) = v_{\mf p_\theta}(b_1) \cdot e_\theta, \quad v(b_2) = v_{\mf p_\eta}(b_2) \cdot e_\eta,
\] 
where $e_\theta$ is the ramification index of $\mf p$ over $\mf p_\theta$, and similar for $e_\eta$. 

We define $K_V$ as in the previous section: recall that $K_V$ is the smallest field extension of $H$ containing all coordinates of all zeros of $\nu$ and $\delta$. Also recall from section \ref{sec:DrinfeldModules} that we have fixed an embedding $\overline K \hookrightarrow \overline{K_v}$. This gives a valuation 
$w$ on $K_V$. We normalize $w$ so that the value group is $\bb Z$. Then once again, for a function $b_3 \in K$, by viewing $b_3 \in \mathrm{Frac} R$, we have 
\[
w(b_3) = v(b_3) \cdot e_w,
\]
where $e_w$ is the ramification index of $w$ over $v$. 

Finally, we let $I_k$ be the $A$-ideal $(\theta - \theta^{q^k}, \eta - \eta^{q^k})$, $J_{k,i}$ be the $R$-ideal $(\theta^{q^k} - \alpha_i, \eta^{q^k} - \beta_i)$, $J_{k,i}'$ be the $R$-ideal $(\theta - \alpha_i'^{q^k}, \eta - \beta_i'^{q^k})$, $J_k$ be the $R$-ideal $\prod_i J_{k,i}$, and $J_k'$ be the $R$-ideal $\prod_i J_{k,i}'$. Hence as fractional $R$-ideals,
\[
(e_n) = \lrp{\prod_{k=0}^{n-1} I_{n-k}^{(k)} R }^{-1} J_n (J_0^{(n)})^{-1},
\]
\[
(l_n) = \lrp{\prod_{k=1}^n I_k R}^{-1} J_{n+1}' J_1'^{-1}.
\]

Thus, 
\begin{prop}
\label{prop:formulaWCoeff}
\[
w(e_n) = w(J_n) - w(J_0^{(n)}) - e_w \sum_{k=0}^{n-1} v(I_{n-k}^{(k)}),
\]
\[
w(l_n) = w(J_{n+1}') - w(J_1') - e_w \sum_{k=1}^n v(I_k).
\]
\end{prop}


\subsection{Main term for the logarithm, and one term for the exponential: \texorpdfstring{$I_k$}{I sub k}}

The main contribution of $w(l_n)$ will come from $v(I_k)$, which we shall first compute. We will also compute $v(I_{n-k}^{(k)})$ for the exponential series. 

\begin{lem}
\label{lem:vForThetaQK}
\[
v_{\mf p_{\theta}} (\theta - \theta^{q^k}) = 
\begin{cases}
1 &\text{ if } \deg_\theta \mf p_\theta \mid k, \\
0 &\text{ else. }
\end{cases}
\]
The same is true for $\eta$. 
\end{lem}

\begin{proof}
All elements in $\bb F_{q^k}$ are roots of the equation $X^{q^k} - X = 0$, so $X^{q^k} - X = \prod_{c \in \bb F_{q^k}} (X-c)$. By grouping $\Gal(\bb F_{q^k}/\bb F_q)$-conjugates together, we have that 
\[
X^{q^k} - X = \prod_{\substack{a \in \bb F_q[X] \\ a \text{ monic irreducible} \\ \deg a \mid k}} a.
\]

\end{proof}

\begin{prop}
\label{prop:vForISubk}
\[
v(I_k) = \begin{cases}
\min\{e_\theta, e_\eta\} &\text{ if } \deg \mf p \mid k \gcd(f_\theta, f_\eta), \\
0 &\text{ else. }
\end{cases}
\]
\end{prop}

\begin{proof}
To have $v(I_k) > 0$, we must have $\mf p$ dividing both $\theta - \theta^{q^k}$ and $\eta - \eta^{q^k}$. By the previous lemma, the first requirement is satisfied when $\deg_\theta \mf p_\theta \mid k$, and the second is satisfied when $\deg_\eta \mf p_\eta \mid k$. Expressing both of these with $\deg \mf p$, we have that $\deg \mf p$ divides both $kf_\theta$ and $kf_\eta$, which is equivalent to $\deg \mf p$ dividing $k\gcd(f_\theta, f_\eta)$. 

Now suppose $v(I_k) > 0$. Recall that 
\[
v(\theta - \theta^{q^k}) = v_\theta(\theta - \theta^{q^k}) \cdot e_\theta,
\]
which is $e_\theta$ by the previous lemma and our assumption that $v(I_k) > 0$. Similarly, $v(\eta - \eta^{q^k}) = e_\eta$. As a result, 
\[
v(I_k) = v((\theta - \theta^{q^k}, \eta - \eta^{q^k})) = \min\{e_\theta, e_\eta\}.
\]

\end{proof}

By adding these up for $k = 1, \cdots, n$, we can see that 
\begin{cor}
\label{cor:sumOfvForISubk}
\[
\sum_{k=1}^n v(I_k) = \min\{e_\theta, e_\eta\} \left\lfloor \frac{n \gcd(f_\theta, f_\eta)}{\deg \mf p} \right\rfloor.
\]
\qed
\end{cor}

We observe in particular that this has order of magnitude $O(n)$. 

We now proceed to compute $v(I_{n-k}^{(k)})$. 

\begin{prop}
For $0 \leq k \leq n-1$,
\[
v(I_{n-k}^{(k)}) = \begin{cases}
\min\{ e_\theta, e_\eta \} q^k &\text{ if } \deg \mf p \mid (n-k)\gcd(f_\theta, f_\eta) \\
0 &\text{ else. }
\end{cases}
\]
\end{prop}

\begin{proof}
By Lemma \ref{lem:vForThetaQK}, 
\[
v_{\mf p_\theta}(\theta^{q^n} - \theta^{q^k})
= v_{\mf p_\theta}((\theta^{q^{n-k}} - \theta)^{q^k})
= \begin{cases}
q^k &\text{ if } \deg_\theta \mf p_\theta \mid k \\
0 &\text{ else. }
\end{cases}
\]
The same is true for $\eta$. The result now follows by a similar argument as Proposition \ref{prop:vForISubk}. 

\end{proof}

For our purpose, we will only give an upper bound for the sum of the terms $v(I_{n-k}^{(n-k)})$. 
\begin{cor}
\label{cor:vForITwisted}
Let $C = \dfrac{\deg \mf p}{\gcd(f_\theta, f_\eta)}$. Then
\[
\sum_{k=0}^{n-1} v(I_{n-k}^{(k)})
< q^n \min\{e_\theta, e_\eta\} \frac{1}{q^C - 1}.
\]
\end{cor}

\begin{proof}
\begin{align*}
\sum_{k=0}^{n-1} v(I_{n-k}^{(k)}) 
&= \sum_{k=1}^{n} v(I_{k}^{(n-k)}) \\
&= \min\{e_\theta, e_\eta\} \lrp{ q^{n - C} + q^{n - 2C} + \cdots + q^{n - \lfloor \frac{n}{C \gcd(f_\theta, f_\eta)} \rfloor C} } \\
&= \min\{e_\theta, e_\eta\} q^n \frac{1 - q^{-C \lfloor \frac{n}{C \gcd(f_\theta, f_\eta} \rfloor}}{q^C - 1} \\
&< q^n \min\{e_\theta, e_\eta\} \frac{1}{q^C - 1}.
\end{align*}

\end{proof}

Even though we give this as an upper bound, we can see from the proof that the sum of the terms $v(I_{n-k}^{(k)})$ has order $O(q^n)$. 


\subsection{Other terms for the logarithm}

We then proceed to show that the term $w(J_n')$ does not matter when compared to the sum of $v(I_k)$'s. 

\begin{prop}
\label{prop:wForJSubniPrime}
Fix an $i$. Then $w(J_{n,i}')$ is bounded by some constant independent of $n$.
\end{prop}

\begin{proof}
If $w(J_{n,i}') > 0$ for only finitely many $n$, we are done. Suppose not. 
Let $N_1 > N_2 > 0$ be two integers where $w(J_{N_1,i}')$ and $w(J_{N_2,i}')$ are positive. Then 
\[
w\lrp{(\theta - \alpha_i'^{q^{N_1}}, \eta - \beta_i'^{q^{N_1}}, \theta - \alpha_i'^{q^{N_2}}, \eta - \beta_i'^{q^{N_2}})} = w(J_{N_1,i}' + J_{N_2,i}') > 0.
\]
We can rewrite this $R$-ideal as
\[
(\theta - \alpha_i'^{q^{N_2}}, \theta - \theta^{q^{N_1 - N_2}}, \eta - \beta_i'^{q^{N_2}}, \eta - \eta^{q^{N_1 - N_2}}),
\]
by considering $\theta - \alpha_i'^{q^{N_1}}$ modulo $\theta - \alpha_i'^{q^{N_2}}$, and similarly for $\eta, \beta_i'$. In particular, we can see that 
\[
J_{N_1,i}' + J_{N_2,i}' \supset I_{N_1-N_2} R.
\]
Thus 
\[
\min\{w(J_{N_1,i}'), w(J_{N_2,i}')\} = w(J_{N_1,i}' + J_{N_2,i}') \leq w(I_{N_1-N_2}).
\]
By Proposition \ref{prop:vForISubk}, the last term is bounded by $e_w \cdot \max\{e_\theta, e_\eta\}$. What we have shown is that if we pick any two distinct positive $N_1, N_2$ with $w(J_{N_1, i}'), w(J_{N_2, i}')$ both positive, then at least one of them is bounded by $e_w \cdot \max\{e_\theta, e_\eta\}$. Thus we can have at most one $N$ such that $w(J_{N,i}')$ is bigger than the $e_w \cdot \max\{e_\theta, e_\eta\}$. Therefore, $w(J_{n,i}')$ is bounded by a constant independent on $n$. 

\end{proof}

\begin{cor}
\label{cor:wForJSubnPrime}
The term $w(J_n')$ is bounded by some constant independent of $n$.
\end{cor}

\begin{proof}
This comes from $\displaystyle w(J_n') = \sum_{i=1}^g  w(J_{n,i}')$. 
\end{proof}


\subsection{Other terms for the exponential}

As for the exponential, the term $w(J_0^{(n)})$ will actually contribute potentially. 
\begin{lem}
\label{lem:wForJSub0iTwisted}
\[
w(J_{0,i}^{(n)}) = q^{n} w(J_{0,i}),
\]
and hence 
\[
w(J_0^{(n)}) = q^{n} w(J_0). 
\]
\end{lem}

\begin{proof}
\begin{align*}
w(J_{0,i}^{(n)})
&= \min\{ w((\theta - \alpha_i)^{(n)}), w((\eta - \beta_i)^{(n)}) \} \\
&= q^{n} \min\{ w(\theta - \alpha_i), w(\eta - \beta_i) \} \\
&= q^{n} w(J_{0,i}). 
\end{align*}

\end{proof}

This means that in the formula for $w(l_n)$ in Proposition \ref{prop:formulaWCoeff}, the term $w(J_{0,i}^{(n)})$ either does not contribute at all, or it has order $O(q^n)$, which is the same as the sum of the terms $v(I_k^{(n-k)})$. 

The term $w(J_n)$ can be ignored when computing the $w$-adic convergence of exponential: since it is positive, it only makes the convergence easier. 

%

%


\subsection{Conclusion for \texorpdfstring{$v$}{v}-adic convergence of exponential and logarithm}

Recall from Proposition \ref{prop:formulaWCoeff} that 
\[
w(l_n) = w(J_{n+1}') -w(J_1') - e_w \sum_{k=1}^n v(I_k).
\]
By using Corollary \ref{cor:sumOfvForISubk} and Corollary \ref{cor:wForJSubnPrime}, we now have the following result. 

\begin{thm}
For $n \gg 0$, $w(l_n)$ is negative, and $|w(l_n)|$ has order of magnitude $O(n)$. \qed
\end{thm}

This allows us to show the $v$ (or equivalently $w$)-adic convergence of Hayes logarithm in $\bb C_v$. 

\begin{thm}
\label{thm:vadicConvOfLog}
The logarithm series for a Hayes module, in the cases of an elliptic curve or a ramifying hyperelliptic curve with $\infty$ the rational point at infinity, converges $w$-adically for all $z \in \bb C_v$ with $w(z) > 0$. 
\end{thm}

\begin{proof}
Suppose $z \in \bb C_v$ with $w(z) > 0$. Then 
\[
w(l_n z^{q^n}) = w(l_n) + q^n w(z).
\]
Since $w(z) > 0$ and $w(l_n)$ has order of magnitude $O(n)$, $w(l_n z^{q^n})$ goes to infinity as $n$ does, showing that the series 
\[
\log_{\rho} (z) = \sum_{n = 0}^\infty l_n z^{q^n}
\]
converges in $\bb C_v$. 
\end{proof}

As for the exponential, we can conclude the following. 
\begin{thm}
\label{thm:vadicConvOfExp}
The exponential series for a Hayes module, in the case of an elliptic curve or a ramifying hyperelliptic curve with $\infty$ as the rational point at infinity, converges $w$-adically for all $z \in \bb C_v$ with 
\[
w(z) > w(J_0) + e_w \min\{e_\theta, e_\eta\} \frac{1}{q^{C}-1},
\]
where $C = \dfrac{\deg \mf p}{\gcd(f_\theta,f_\eta)}$. 
\end{thm}

\begin{proof}
By combining Corollary \ref{cor:vForITwisted}, and Lemma \ref{lem:wForJSub0iTwisted}, we can see that
\begin{align*}
w(e_n) 
&= w(J_n) - w(J_0^{(n)}) - e_w \sum_{k=0}^{n-1} v(I_{n-k}^{(k)}) \\
&\geq - \lrp{w(J_0) + e_w \min\{e_\theta, e_\eta\} \frac{1}{q^C-1} } q^n. 
\end{align*}
Hence if $w(z)$ is strictly bigger than the number in the parentheses, we have that
\[
w(e_n z^{q^n}) \geq \e q^n,
\]
for some $\e > 0$. This valuation goes to infinity as $n \to \infty$. This shows that the series 
\[
e_{\rho}(z) = \sum_{n=0}^\infty e_n z^{q^n}
\]
converges in $\bb C_v$. 
\end{proof}

\begin{defn}
We denote by $e_{v,\rho}(z)$ and $\log_{v,\rho}(z)$ the functions defined by the respective series on $\bb C_v$, whenever the series converges. 
\end{defn}

%% file: Sec6.tex
As an application of the $v$-adic convergence of $\log_{v,\rho}$, we prove a log-algebraicity formula for the $v$-adic Goss $L$-value $L_v(1,\Psi)$. This follows the philosophy of \cite[(43)]{Anderson1996} in Anderson's paper and uses Lutes' result \cite{Lutes2010} on log-algebraicity for $L(1,\Psi)$. 

\subsection{Goss \texorpdfstring{$L$}{L}-function}
We first define the $L$-functions we are interested in. These characteristic $p$ $L$-functions for function fields were first studied by Carlitz \cite{Carlitz1935} (for zeta functions on Carlitz module), then Goss \cite{Goss1978_vonStaudt}, \cite{Goss1979}, \cite{Goss1980_piAdicEisSeries}, etc (for zeta / $L$-functions and $v$-adic zeta / $L$-functions on a general curve). We will only define the functions with domain the set of integers. For a much more general definition, see \cite[\S8]{Goss1998_BasicStrcBook}. 

Let us begin with the assumption that $A$ is a PID, and consider just the zeta function instead of $L$-functions. In this case, the definition of the Goss zeta function is analogous to that of the Riemann zeta function. We use $A^+$ to denote the subset of elements of $\sgn$ 1 in $A$, where ``$\sgn$'' is taken with respect to ``$\infty$'' in $K$, not $\bb K$. We also use $A_d^+$ to denote the subset of $\sgn$ 1 elements of degree $d$.

\begin{defn}
For $A$ a PID, the {\em Goss zeta function} is 
\[
\zeta(n) = \sum_{a \in A^+} \frac{1}{a^n}
:= \sum_{d \geq 0} \sum_{a \in A_d^+} \frac{1}{a^n}.
\]
The second summation indicates the order in which the terms should be summed. 
\end{defn}

For $A$ not a PID, in order to define the zeta function, we need to define the exponentiation of an ideal by an integer. The definition is given carefully in Goss's book \cite[\S8.2]{Goss1998_BasicStrcBook}, and we refer the readers to this for more details.

\begin{defn}
For an ideal $I \subset A$ and a positive integer $n$, we denote by $I^{[n]}$ an element in a finite (not necessarily separable) extension of $K$, satisfying the properties given in Goss \cite[\S8.2]{Goss1998_BasicStrcBook}. In particular, if $I$ is a principal ideal generated by a $\sgn$ 1 element $i$, then $I^{[n]} = i^n$. 
\end{defn}

From \cite[8.2.8]{Goss1998_BasicStrcBook}, we see that the finite (not necessarily separable) extension of $K$ in the definition above can be chosen to not depend on $I$ or $n$. We denote that extension by $K'$. 

This gives a natural extension of zeta functions to $A$ a non-PID. Let $\mc I$ be the group of fractional $A$-ideals in $K$, $\mc P^+$ the subgroup of principal ideals generated by $\sgn$-1 elements, $h^+ = \# \mc I / \mc P^+$ the narrow class number, and $\mf a_1, \cdots, \mf a_{h^+}$ representatives of the narrow class group $\mc I / \mc P^+$. We borrow the notation from number fields and use $\mf a_j^{-1}$ to stand for the inverse fractional $A$-ideal of $\mf a_j$. We use $\mf a_j^{-1,+}$ to denote the subset of $\sgn 1$ elements in $\mf a_j^{-1}$. 
\begin{defn}
The {\em Goss zeta function} is defined as 
\[
\zeta(n) = \sum_{I \subset A} \frac{1}{I^{[n]}} := \sum_{j=1}^{h^+} \frac{1}{\mf a_j^{[n]}} \sum_{d = 0}^\infty \sum_{\substack{i \in \mf a_j^{-1,+} \\ \deg i + \deg \mf a_j = d}} \frac{1}{i^n}.
\]
Once again, the second summation indicates how the terms should be summed. 
\end{defn}

As for the {\em $v$-adic zeta function}, it is defined in a very similar manner, but we discard those $i$ such that $i \mf a_j$ is divisible by $\mf p$. 

Finally, if $\Psi$ is a character on $\mc I$ of conductor $\mf p$, we can define the {\em Goss $L$-function} and its $v$-adic analogue very similarly:
\begin{defn}
For $n \in \bb Z$, the {\em Goss $L$-function} is defined as
\[
L(n,\Psi) = \sum_{0 \neq I \subset A} \frac{\Psi(I)}{I^{[n]}} := \sum_{j = 1}^{h^+} \frac{1}{\mf a_j^{[n]}} \sum_{d \geq 0} \sum_{\substack{i \in \mf a_j^{-1,+} \\ \deg i + \deg \mf a_j = d}} \frac{\Psi(i\mf a_j)}{i^n}.
\]
The {\em $v$-adic Goss $L$-function} $L_v(n,\Psi)$ is defined by discarding the terms with $v(i\mf a_j) > 0$. 
\end{defn}

It is a remarkable result by Goss that all such series converge for all integers $n$, in $\bb C_\infty$ and in $\bb C_v$ respectively. The proof uses a lemma \cite[8.8.1]{Goss1998_BasicStrcBook}, which gives a $v$-adic estimation for the degree $d$ partial sum. 


\subsection{Anderson's log-algebracity}

From now on we focus on the values $L(1,\Psi)$ and $L_v(1,\Psi)$. Classically over $\bb Q$ with $\chi$ an even non-trivial character, both analogues of $L$-series have an expression involving logs of algebraic numbers. That is, we have that 
\[
L(1,\chi) = \sum \alpha_i \log S_i, \qquad 
L_p(1,\chi) = \sum \alpha_i' \log_p S_i'
\]
for some algebraic numbers $\alpha_i, S_i, \alpha_i', S_i' \in \overline{\bb Q}$. See \cite[4.9, 5.18]{Washington1997}. 

As for function fields, Anderson has proved similar log-algebraic formulas for both $L$-functions and $v$-adic $L$-functions on $\bb F_q[\theta]$, and Lutes proved a formula for $L$-functions on all curves. We shall first briefly go through their results. For a general curve $X$, let $H^+$ be extension of $H$ that is the field of definition for a Hayes module $\rho$, see \cite[\S3]{Thakur2004}. Let $\rho_I(\tau)$ be the monic generator of the principal left ideal $\{\rho_i(\tau) ~|~ i \in I\} \subset H^+\{\tau\}$, $D(\rho_I)$ be its constant term, and $\rho_I(Y)$ be obtained by replacing $\tau^j$ with $Y^{q^j}$ in $\rho_I(\tau)$. Here $Y$ is a free variable, in particular transcendental over $H^+$. We denote by $\rho[I]$ the $I$-torsion points of $\rho$. 

Let $\Lambda_\rho$ be the lattice corresponding to $\rho$, $\varpi_\rho$ a generator of the lattice as an $A$-module. To ease notations, we set $\mathbf{e}_A(z) := e_\rho(\varpi_\rho z)$, and $\mathbf{e}_I(z) := e_{I * \rho}(D(\rho_I) \varpi_\rho z)$.

\begin{defn}
\begin{enumerate}
\item Let 
\[
b(Y) = \sum_{i = 0}^{\deg b} b_i Y^i \in H^+[Y]. 
\]
We define an action of nonzero $A$-ideals on $H^+[Y]$ by
\[
(J * b)(Y) := \sum_{i=0}^{\deg b} b_i^{(J, H^+/K)} \rho_J(Y)^i,
\]
where $(J,H^+/K) \in \Gal(H^+/K)$ is the Artin symbol. 

\item 
For each $b \in H^+[Y]$, let $l(b;z)$ be a power series in $z$ over $H^+[Y]$ defined by 
\[
l(b; z) := \sum_J \frac{J * b}{D(\rho_J)} z^{q^{\deg J}},
\]
where the sum is over all nonzero ideals $J \subset A$. Similarly, upon fixing a nonzero ideal $I \subset A$, we define 
\[
l_I(b;z) := \frac{1}{D(\rho_I)} \sum_{\alpha \in (I^{-1,+}} \frac{(\alpha I) * b}{\alpha} z^{q^{\deg I + \deg \alpha}}.
\]
\end{enumerate}
\end{defn}

As a remark, it is clear that $l(b;z) = {\displaystyle \sum_{i=1}^{h^+} l_{\mf a_i}(b;z)}$, where $\mf a_i$ goes over all classes of $\mc I / \mc P^+$. Anderson's main result in \cite{Anderson1996} asserts that given $b \in \mc O_{H^+}[Y]$, $l(b;z)$ is some sort of logarithm. To be precise, we apply the exponential power series associated to $\rho$, and can obtain a polynomial over $\mc O_{H^+}$. 

\begin{thm}
\label{AndersonLogAlg}
(Anderson, \cite[Theorem 3]{Anderson1996}) For $b \in \mc O_{H^+}[Y]$, the formal power series 
\[
S(b;z) := e_\rho (l(b;z)),
\]
a priori in the ring $(H^+[Y])[[z]]$, is in fact in $\mc O_{H^+}[Y,z]$. 
\end{thm}

This is an analogue of the fact that $\exp(\log(1-z))$, a priori a formal power series in $\bb Q[[z]]$, is in fact in $\bb Z[z]$. Readers can refer to \cite{Anderson1996} for a more careful formulation of this theorem, via writing $l(b;z)$ using increasing powers of $z$ and defining $S(b;z)$ in terms of such coefficients. Readers can also refer to a beautiful survey article by Perkins \cite{Perkins13}. 

By definition of $S(b;z)$, we have that $S(Y^m;z)$ is divisible by $Y^m$. We will use this observation very soon.





\subsection{\texorpdfstring{$L(1, \chi)$}{L(1, chi)}}

Next, we pick appropriate $b \in \mc O_{H^+}[Y]$, and evaluate at particular values of $Y$ and $z$. 

\begin{defn}
\label{defnOflmI}
For a nonzero ideal $I \subset A$, we define 
\[
l_{m,I}(Y) = \frac{1}{D(\rho_I)} \sum_{\alpha \in I^{-1,+}} \frac{\mathbf{e}_I(\alpha x)}{\alpha}
= l_I(Y^m; z)|_{Y = \mathbf{e}_A(x), z = 1}.
\]
Here $x$ is another formal variable. Later we will substitute $x$ by elements in $\mf p^{-1} \Lambda_\rho / \Lambda_\rho$.
\end{defn}

Fix an isomorphism 
\[
A / \mf p \longrightarrow \mf p^{-1} \Lambda_\rho / \Lambda_\rho
\]
\[
a \mapsto a\mu
\]

Recall that $\mf a_1, \dots, \mf a_{h^+} \subset A$ are ideals representing the classes of $\mc I / \mc P^+$. We then substitute $b = Y^m$, $Y = \mathbf e_A(a\mu)$, and $z = 1$ in Anderson's log-algebraicity theorem and obtain that 
\[
e_\rho\left( \sum_{i = 1}^{h^+} l_{m,\mf a_i}(a\mu) \right) = S(Y^m;z)|_{Y = \mathbf e_A(a\mu), z = 1}
\]
is an algebraic integer in $\mc O_{H^+}[\rho[\mf p]]$. Before we move on to values of $L(1,\chi)$, we shall first mention an important lemma. This lemma is a key to prove the log-algebraicity for $L_v(1,\chi)$. 

\begin{lem}
\label{SYmzDivisibleBymPower}
For any $\beta \in \mc O_{K(\rho[\mf p])}$, $S(\beta Y^m;z)|_{Y = \mathbf e_A(a\mu), z = 1}$ is divisible by $\mathbf e_A(a\mu)^m$ in $\mc O_{K(\rho[\mf p])}$. 
\end{lem}

\begin{proof}
Recall from the discussion after Theorem \ref{AndersonLogAlg} that $S(Y^m; z)$ is divisible by $Y^m$. The exact same argument shows that the same is true for $S(\beta Y^m;z)$. Now evaluate at $Y = \mathbf e_A(a\mu), z = 1$, we have that $S(\beta Y^m;z)|_{Y = \mathbf e_A(a\mu), z = 1}$ is divisible by 
\[ \mathbf e_A(a\mu)^m = (\rho_a(\mathbf e_A(\mu)))^m, \]
which is divisible by $\mathbf e_A(\mu)^m$.
\end{proof}

We now return to the study of $L(1,\Psi)$. With the notation $l_{m,I}(x)$, Anderson and Lutes give formula for $L(1,\Psi)$ in terms of $\mathbf e_I(a\mu)$'s and $l_{m,I}(b\mu)$, using Lagrange interpolation. 

\begin{thm}
\label{AndersonLutesL1Chi}
\begin{enumerate}
\item (Anderson, \cite[(38)]{Anderson1996})
Let $A = \bb F_q[\theta]$ and $\chi: A \to \bb C_\infty$ a character of conductor $\mf p$. In this case $\mc I / \mc P^+ = 1$. Then 
\[
L(1,\chi) = -\frac{1}{\mf p^{[1]}} \sum_{m=1}^{q^{\deg \mf p}-1} \left( \sum_{a \in \bb F_{\mf p}^\times} \mathbf e_{m}^*(a) \right) \left( \sum_{b \in \bb F_{\mf p}^\times} \chi^{-1}(b) l_{m,A}(b\mu) \right).
\]
The $\mathbf e_m^*(a)$ are algebraic, obtained via Lagrange interpolation. 

\item (Lutes, \cite[V.13]{Lutes2010})
For a general $A$, and $\Psi$ a character of conductor $\mf p$ on $\mc I(\mf p)$, the group of fractional $A$-ideals prime to $\mf p$, 
\begin{multline*}
L(1,\Psi) 
= \sum_{j = 1}^{h^+} {\Bigg[} -\frac{\Psi(\mf a_j)}{\mf a_j^{[1]}} \sum_{m=1}^{q^{\deg \mf p}-1}
\sum_{a \in \bb F_{\mf p}^\times} \left( \Psi(a) \mathbf e_{m,\mf a_j}^*(a/\nu_j) \right) \\
\sum_{b \in \bb F_{\mf p}^\times} \left( \Psi^{-1}(b) l_{m,\mf a_j}(b\mu) \right)
{\Bigg]}.
\end{multline*}

The $\nu_j$'s, $\mathbf e_{m,\mf a_j}(a/\nu_j)$'s are defined similarly as $\mu$ and $\mathbf e_m^*(a)$. See \cite{Lutes2010} for detailed definitions. 

\end{enumerate}
\end{thm}

When $A = \bb F_q[\theta]$, we can rewrite the term $l_{m,A}(b\mu)$ as 
\[
l_{m,A}(b\mu) = l(Y^m;z)|_{Y = \mathbf{e}_A(b\mu), z = 1} = \log_\rho S(Y^m;z)|_{Y = \mathbf{e}_A(b\mu), z = 1},
\]
in particular the Carlitz logarithm of an algebraic integer. This shows that $L(1,\chi)$ is log-algebraic. 

However, in Lutes's scenario things are not so immediate. A direct difficulty we encounter is that we cannot immediately see that $l_{m,\mf a_j}(b\mu)$ is log-algebraic. To fix this, Lutes fixes a set of $K$-linearly independent elements $\{\beta_1, \cdots, \beta_{h^+}\}$ in $\mc O_{H^+}$. Then he applies Anderson's log-algebraicity theorem \ref{AndersonLogAlg} to $b = \beta_j Y^m$ instead of $Y^m$. This tells us that 
\[
S(\beta_j Y^m;z) = e_\rho(l(\beta_j Y^m ;z))
\]
is a polynomial in $Y, z$ with coefficients in $\mc O_{H^+}$. Evaluating this at $Y = \mathbf e_A(b \mu), z = 1$, we have that
\[
S(\beta_j Y^m;z)_{Y = \mathbf e_A(b \mu), z = 1}
\]
is an algebraic integer for all $\beta_j$. Now, unwinding the right hand side, we see that
\begin{align*}
S(\beta_j Y^m;z)_{Y = \mathbf e_A(b \mu), z = 1}
&= e_\rho \left( l(\beta_j Y^m ; z) \right)_{Y = \mathbf e_A(b \mu), z = 1} \\
&= \sum_{i=1}^{h^+} e_\rho \left( l_{\mf a_i}( \beta_j Y^m; z) \right)_{Y = \mathbf e_A(b \mu), z = 1} \\
&= \sum_{i=1}^{h^+} e_\rho \left( \beta_j^{(\mf a_i, H^+ / K)} l_{\mf a_i}(Y^m;z) \right)_{Y = \mathbf e_A(b \mu), z = 1} \\
&= \sum_{i=1}^{h^+} e_\rho( \beta_j^{(\mf a_i, H^+ / K)} l_{m,\mf a_i}(b\mu) ). 
\end{align*}
As $\beta_j$ varies in the set $\{\beta_1,\dots,\beta_{h^+}\}$, these equations can be written in terms of a matrix equation. For ease of notation, set $\mf L_j := \log_\rho \lrp{S(\beta_j Y^m;z)_{Y = \mathbf e_A(b \mu), z = 1}}$. 
\[
\begin{pmatrix}
\mf L_1 \\ \mf L_2 \\ \vdots \\ \mf L_{h^+}
\end{pmatrix}
= \begin{pmatrix}
\beta_1^{(\mf a_1, H^+/K)} & \beta_1^{(\mf a_2, H^+/K)} & \cdots & \beta_1^{(\mf a_{h^+}, H^+/K)} \\ 
\beta_2^{(\mf a_1, H^+/K)} & \beta_2^{(\mf a_2, H^+/K)} & \cdots & \beta_2^{(\mf a_{h^+}, H^+/K)} \\
\vdots & \vdots & \ddots & \vdots \\  
\beta_{h^+}^{(\mf a_1, H^+/K)} & \beta_{h^+}^{(\mf a_2, H^+/K)} & \cdots & \beta_{h^+}^{(\mf a_{h^+}, H^+/K)}
\end{pmatrix}
\begin{pmatrix}
l_{m,\mf a_1}(b\mu) \\ l_{m,\mf a_2}(b\mu) \\ \vdots \\ l_{m,\mf a_{h^+}}(b\mu)
\end{pmatrix}.
\]
By definition, $\{\beta_1, \cdots, \beta_{h^+}\}$ is $K$-linearly independent. Thus the matrix here has linearly-independent rows, equivalently nonzero determinant. This shows that $l_{m,\mf a_i}(b\mu)$ is an $H^+$-linear combination of $\mf L_j = \log_\rho \left( S(\beta_j Y^m;z)_{Y = \mathbf e_A(b\mu), z=1} \right)$. Combining this with theorem \ref{AndersonLutesL1Chi}, we obtain the following. 

\begin{thm}
\cite[V.14]{Lutes2010}
\label{thm:L(1,Psi)LogAlg}
$L(1,\Psi)$ is log-algebraic, i.e. there exists $\alpha_1, \cdots, \alpha_s$, $S_1, \cdots, S_s \in \overline K$ such that 
\[
L(1,\Psi) = \sum_i \alpha_i \log_\rho S_i.
\]
In fact, we can pick $S_i \in K(\rho[\mf p])$, $\alpha_i \in K'(\rho[\mf p])$, and $s \leq (q^{\deg \mf p}-1)^2 \cdot h^+$. \qed
\end{thm}

If one wishes, the $\alpha_i$'s can be written very explicitly, by using Cramer's rule to solve the matrix equation. The expression will get too long, so we skip it here. 











\subsection{Log-algebraicity for \texorpdfstring{$L_v(1,\Psi)$}{v-adic L(1, Psi)} on Elliptic or Ramifying Hyperelliptic curves}

We now proceed to prove our application, which is the log-algebraicity for $L_v(1,\Psi)$ for our curves. The following proposition will be imitating \cite[Proposition 12]{Anderson1996} by Anderson. The idea is that we can use the $v$-adic convergence of $\log_{v,\rho}$ to make a formal series converge to what we desire. 

\begin{prop}
\label{prop:vadicConvForlmI}
Fix $\beta \in \mc O_{H^+}$. The series 
\begin{multline*}
\sum_{j=1}^{h^+} \beta^{(\mf a_j, H^+/K)} l_{m,\mf a_j}(b\mu) \\
= \sum_{k = 0}^\infty \left( \sum_{j = 1}^{h^+} \frac{1}{D(\rho_{\mf a_j})} \beta^{(\mf a_j, H^+/K)} \left( \sum_{\substack{\alpha \in (\mf a_j^{-1,+} \\ \deg \mf a_j + \deg \alpha = k}} \frac{\mathbf e_{\mf a_j}(\alpha b\mu)^m}{\alpha} \right) \right),
\end{multline*}
summing in the order as indicated, converges $v$-adically to 
\[
\log_{v,\rho} {\big(}S(\beta Y^m;z)|_{Y = \mathbf e_A(b\mu), z = 1} {\big)}. 
\]
\end{prop}

\begin{proof}
By definition of $e_I(z)$, 
\[
\mathbf e_{\mf a_j}(\alpha b \mu)
= e_{\mf a_j * \rho}(D(\rho_{\mf a_j}) \varpi_\rho \alpha b\mu). 
\]
We shall investigate the right hand side. Since $\alpha$ is positive, by definition of $\rho_I$ we have that 
\[
D(\rho_{\alpha \mf a_j}) = \alpha D(\rho_{\mf a_j}).
\]
Since $\alpha \mf a_j$ is in the same ideal class as $\mf a_j$, $(\alpha \mf a_j)*\rho$ and $\mf a_j * \rho$ are the same Drinfeld module. As a result
\[
e_{(\alpha \mf a_j)*\rho} = e_{\mf a_j * \rho}. 
\]
Thus 
\[
\mathbf e_{\mf a_j}(\alpha b \mu)
= e_{(\alpha \mf a_j) * \rho}(D(\rho_{\alpha \mf a_j}) \varpi_\rho b\mu)
= \rho_{\alpha \mf a_j} (\mathbf e_A(b\mu)),
\]
where functional equation for $\rho_I$ is used in the last equality. Putting everything together, the degree-$k$ part of our series becomes
\begin{align*}
&\phantom{==} \sum_{j = 1}^{h^+} \sum_{\substack{\alpha \in (\mf a_j^{-1,+} \\ \deg \mf a_j + \deg \alpha = k}} \frac{\beta^{(\mf a_j, H^+/K)}}{\alpha D(\rho_{\mf a_j})} (\rho_{\alpha \mf a_j}(\mathbf e_A(b\mu)))^m \\
&= \left. \sum_{j = 1}^{h^+} \sum_{\substack{\alpha \in (\mf a_j^{-1,+} \\ \deg \mf a_j + \deg \alpha = k}} \frac{1}{D(\rho_{\alpha \mf a_j})}  ((\alpha \mf a_j) * (\beta Y^m)) \right|_{Y = \mathbf e_A(b\mu)}.
\end{align*}
The sum goes over all ideals $J \subset A$ with $\deg J = k$. Hence this equals
\[
\left. \sum_{\deg J = k} \frac{J * (\beta Y^m)}{D(\rho_J)} \right|_{Y = \mathbf e_A(b\mu)}.
\]
Therefore, our whole series is the same as 
\[
\sum_{k = 0}^\infty \left. \left( \sum_{\deg J = k} \frac{J * (\beta Y^m)}{D(\rho_J)} \right) \right|_{Y = \mathbf e_A(b\mu)}
= l(\beta Y^m; z)|_{Y = \mathbf e_A(b\mu), z = 1}.
\]
By Anderson's log-algebraicity theorem \cite{Anderson1996}, the formal power series 
\[
e_\rho(l(\beta Y^m; z))
\]
is actually a polynomial in $Y,z$ over $\mc O_{H^+}$. Hence as formal series, 
\begin{align*}
l(\beta Y^m; z)|_{Y = \mathbf e_A(b\mu), z = 1} 
&= \log_{v,\rho} S(\beta Y^m;z)|_{Y = \mathbf e_A(b\mu), z = 1} \\
&= \sum_{n \geq 0} L_n (S(\beta Y^m;z)|_{Y = \mathbf e_A(b\mu), z = 1})^{q^n},
\end{align*}

\noindent where $L_n$ are the coefficients of logarithm. To makes sense of this series, we need the right hand side to converge $v$-adically. By Lemma \ref{SYmzDivisibleBymPower}, $S(\beta Y^m;z)|_{Y = \mathbf e_A(b\mu), z = 1}$ is an algebraic integer divisible by $\mathbf e_A(\mu)^m$. Since $\mathbf e_A(\mu)$ has positive $v$-adic valuation, so is $S(\beta Y^m;z)|_{Y = \mathbf e_A(b\mu), z = 1}$. Therefore, by Theorem \ref{thm:vadicConvOfLog} the series converges $v$-adically to $\log_{v,\rho} (S(\beta Y^m; z)|_{Y = \mathbf e_A(b\mu) z = 1})$. 

\end{proof}

Following the method of Anderson \cite{Anderson1996}, we are now able to prove the log-algebraicity of $L(1,\Psi)$ for elliptic curves and ramifying hyperelliptic curves. 

\begin{thm}
\label{thm:Lv(1,Psi)LogAlg}
$L_v(1,\Psi)$ is log-algebraic, i.e. there exists 
\[ \alpha_1, \cdots, \alpha_s, S_1, \cdots, S_s \in \overline K, \]
with $v(S_i) > 0$, such that 
\[
L_v(1,\Psi) = \sum_i \alpha_i \log_{v,\rho} S_i.
\]
The $\alpha_i$, $S_i$, $s$ are the same as in log-algebraic theorem \ref{thm:L(1,Psi)LogAlg} for $L(1,\Psi)$. 
\end{thm}

\begin{proof}
The expression 
\[
\sum_{j=1}^{h^+} \beta^{(\mf a_j, H^+/K)} l_{m,\mf a_j}(b\mu)
\]
we considered in Proposition \ref{prop:vadicConvForlmI} is precisely $S(\beta Y^m; z)_{Y = \mathbf e_A(b\mu), z = 1}$. As in the proof of Theorem \ref{thm:L(1,Psi)LogAlg}, let $\{\beta_1, \cdots, \beta_{h^+}\}$ be a $K$-linearly independent subset of $\mc O_{K(\rho[\mf p])}$. By going through the same argument as the proof of Theorem \ref{thm:L(1,Psi)LogAlg} again, but doing everything $v$-adically, we arrive at the desired equality. 
\end{proof}


%% file: Sec7.tex

\subsection{\texorpdfstring{$\bb A = \bb F_3[t,y]/(y^2-(t^3-t-1)), g = 1, h = 1$}{g = 1, h = 1, F3}}
\,\\
(\cite[11.5]{Hayes1979}, \cite[2.3c]{Thakur1993_ShtukasAJ}, \cite[9.1]{GreenPapanikolas2018}, \cite[VIII.4]{Lutes2010}) Our first example is an elliptic curve over $\bb F_3$ with $h(\bb A) = 1$. Let $\bb A = \bb F_3[t,y]/(y^2-(t^3-t-1))$. Then $V = (\theta+1,\eta)$, and 
\[
f = \frac{y - \eta - \eta(t-\theta)}{t - (\theta+1)}.
\]
The Hayes module is given by 
\[
\rho_t = \theta + \eta(\theta^3-\theta) \tau + \tau^2, \, \,
\rho_y = \eta + \eta(\eta^3-\eta) \tau + (\eta^9+\eta^3+\eta) \tau^2 + \tau^3. 
\]
Let $v$ be place on $K$ corresponding to the prime ideal $\mf p := (\theta)$. Fix $\sqrt{-1} \in \bb F_9$. This gives a character of conductor $\mf p$ by $\chi: A \to \bb F_9$, $a(\theta,\eta) \mapsto a(0,\sqrt{-1})$. Let $\lambda \in K(\rho[\mf p])$ be a primitive $t$-torsion point of $\rho$, i.e. a generator of $\rho[\mf p]$ as an $\bb A$-module, and let $\lambda' = \rho_y(\lambda)$. In \cite[VIII.4]{Lutes2010}, a log-algebraic formula of $L(1,\chi)$ is given as 
\begin{align*}
L(1,\chi) &= \frac{\log_\rho(\lambda') + \sqrt{-1} \log_\rho(\lambda)}{\lambda' + \sqrt{-1}\lambda}, \\
L(1,\chi^3) &= \frac{\log_\rho(\lambda') - \sqrt{-1} \log_\rho(\lambda)}{\lambda' - \sqrt{-1}\lambda}. 
\end{align*}

By Theorem \ref{thm:Lv(1,Psi)LogAlg}, we obtain a log-algebraic formula for $L_v(1,\chi)$ given by the same numbers. 
\begin{align*}
L_v(1,\chi) &= \frac{\log_{v,\rho}(\lambda') + \sqrt{-1} \log_{v,\rho}(\lambda)}{\lambda' + \sqrt{-1}\lambda}, \\
L_v(1,\chi^3) &= \frac{\log_{v,\rho}(\lambda') - \sqrt{-1} \log_{v,\rho}(\lambda)}{\lambda' - \sqrt{-1}\lambda}. 
\end{align*}


\subsection{\texorpdfstring{$\bb A = \bb F_2[t,y] / (y^2+y + (t^5+t^3+1), g = 2, h = 1$}{g = 2, h = 1, F2}}
\,\\
(\cite[11.6]{Hayes1979}, \cite[2.3d]{Thakur1993_ShtukasAJ}) This is the only genus at least 2 example with $h(\bb A) = 1$ (cf. \cite{LMQ1975,Stripe2014}). Let $\bb A = \bb F_2[t,y]/(y^2+y+(t^5+t^3+1))$. Then $V = (\theta, \eta+1) + (\theta^2+1, \eta^2+\theta^4)$, and 
\[
f = \frac{y+\eta + (t+\theta)(\theta^4+\theta^3+\theta^2(t+1))}{t^2 + (\theta^2+\theta+1)t + (\theta^3+\theta)}.
\]
The Hayes module is given by 
\[
\rho_t = \theta + (\theta^2+\theta)^2\tau + \tau^2, \qquad
\rho_y = \eta + y_1\tau + y_2\tau^2 + y_3\tau^3 + y_4\tau^4 + \tau^5,
\]
where 
\begin{align*}
y_1 &= (\theta^2+\theta)(\eta^2+\eta) \\
y_2 &= \theta^2(\theta+1)(\eta^2+\eta)(\eta+\theta^3)(\eta+\theta^3+1) \\
y_3 &= \eta(\eta+1)(\theta^5+\theta^3+\theta^2+\theta+1) [(\theta^3+\theta^2+1)\eta + \theta^7+\theta^4+\theta^2] \\
&\phantom{==} [(\theta^3+\theta^2+1)\eta + \theta^7+\theta^4+\theta^3 + 1] \\
y_4 &= [\theta(\eta^2+\eta)(\theta^5+\theta^2+1)(\eta+\theta) (\eta+\theta+1)]^2.
\end{align*}


To illustrate our Corollary \ref{cor:factorizationModXi} and \ref{cor:sumOfvForISubk}, we have factorized the first few coefficients of $\log_\rho$ as $A$-ideals. 
\begin{align*}
(l_1) &= (\theta)(\theta+1), \\
(l_2) &= (\theta)^{-1} (\theta+1)^{-1} (\theta^8+\theta^6+\theta^5+\theta^4+\theta^3+\theta+1), \\
(l_3) &= (\theta^2+\theta+1)^{2}(\theta^{10}+\theta^9+\theta^8+\theta^3+\theta^2+\theta+1), \\
(l_4) &= (\theta^2 + \theta + 1)^{-1} (\theta)^{-2} (\theta+1)^{-2} (\eta + \theta^2) (\eta + \theta^2 + 1) (\eta + \theta) (\eta + \theta + 1) \\
& (\eta + \theta^2 + \theta) (\eta + \theta^2 + \theta + 1) (\eta+1) (\eta) \\
& (\theta^{12} + \theta^9 + \theta^8 + \theta^6 + \theta^3 + \theta^2 + 1) (\theta^4 + \theta + 1).
\end{align*}

As an illustration of Proposition \ref{prop:formulaWCoeff} and Corollary \ref{cor:sumOfvForISubk}, the primes $(\theta)$ and $(\theta+1)$, both of degree 2, divide $l_2^{-1}$ once and $l_4^{-1}$ twice. It is worth noticing that Corollary \ref{cor:sumOfvForISubk} predicts that $(\theta)$ and $(\theta+1)$ should also divide $l_3^{-1}$ once, but they did not show up in the above factorization. This is because both of them get canceled by the terms coming from $V'^{(4)}$, equivalently $w(J_{n+1}')$ as in Proposition \ref{prop:formulaWCoeff}.


\subsection{\texorpdfstring{$\bb A = \bb F_3[t,y] / (y^2 - (t^3+t^2+t), g = 1, h = 2$}{g = 1, h = 2, F3}}
\,\\
(\cite[11.7]{Hayes1979}, \cite[9.2]{GreenPapanikolas2018}, \cite[VIII.5]{Lutes2010}) The third example is a class number 2 elliptic curve over $\bb F_3$. Let $\mathbf A = \bb F_3[t,y] / (y^2 - (t^3 - t^2 - t))$. We have $\mathrm{cl}(\mathbf A) = 2$, $H = K(\sqrt{\theta})$ and $\mc O_H = \bb F_3[\sqrt{\theta}, \frac{\eta}{\sqrt{\theta}}]$. We have to fix a Hayes module with respect to $\sgn$. Set 
\[
\rho_t = \theta + (\sqrt{\theta} - \theta^{\frac{9}{2}} - \eta - \eta^3) \tau + \tau^2,
\]
and $\rho_y$ is determined uniquely from $\rho_t$ by $\rho_t \rho_y = \rho_y \rho_t$. For this Hayes module, we have $V = (- \theta - 1 - \frac{\eta}{\sqrt{\theta}}, -\eta - \theta^{\frac{3}{2}} - \sqrt{\theta})$, and
\[
f = \frac{y - \eta - (-\eta - \theta^{\frac{3}{2}} + \sqrt{\theta})(t - \theta)}{ t + \theta + 1 + \frac{\eta}{\sqrt{\theta}} }.
\]

We have that
\begin{align*}
(l_1^{-1}) = & (\sqrt{\theta}) (\frac{\eta}{\sqrt \theta} - 1, \theta + 1)^{-1}, \\
(l_2^{-1}) = & (\sqrt \theta)^2 (\frac{\eta}{\sqrt \theta} + \sqrt \theta) (\frac{\eta}{\sqrt \theta} - \sqrt \theta)
(\frac{\eta}{\sqrt \theta} + 1, \theta + 1) (\frac{\eta}{\sqrt \theta} - 1, \theta + 1) \\
& (\sqrt \theta + 1) (\sqrt \theta - 1) (\frac{\eta}{\sqrt \theta}) \left( \theta^2(\theta - 1)^3 \frac{\eta}{\sqrt \theta} + (\theta^5 + \theta^3 + \theta^2 + 1) \right)^{-1}
\end{align*}

This time, the ideal coming from $\Xi^{(1)}$ is 
\[
(\theta^2 - \theta, \eta^2 - \eta) = (\sqrt \theta),
\]
instead of $(1)$, since $(\sqrt \theta) \cap A = (\theta)$ in $A$ is of degree 2, not 1. As we can see, $(\sqrt \theta)$ divides $(L_2^{-1})$ twice, matching the prediction from proposition \ref{prop:vForISubk}. 

Let $v$ be the place on $K$ corresponding to the prime ideal $\mf p := (\theta, \eta)$. Then (\cite[11.7]{Hayes1979})
\[
\rho_{\mf p} = \lrp{1 + \theta + \frac{\eta}{\sqrt{\theta}} }^{-1} \sqrt{\theta} + \tau.
\]
We have a character of conductor $\mf p$ given by $\chi: A \to \bb F_3$ by $a(\theta,\eta) \mapsto a(0,0)$. Extend $\chi$ to a character $\Psi$ on the group of fractional ideals in $K$. Fix $\lambda$ to be a primitive $\mf p$-torsion point. A log-algebraic expression for $L(1,\Psi)$ is given by (\cite[VIII.5]{Lutes2010})
\begin{multline*}
L(1,\Psi) = \lrp{\frac{D(\rho_{\mf p})}{\sqrt{\theta}} - \frac{1}{\lambda}} \log_\rho S(X;1)_{X = \lambda} \\
+ \lrp{-\frac{\sqrt{D(\rho_{\mf p}}}{\theta} - \frac{1}{\lambda \sqrt{\theta}} } \log_\rho S(\sqrt{\theta} X; 1)_{X = \lambda}.
\end{multline*}
Lutes also computed the special polynomial $S(X;z)$: 
\begin{multline*}
S(X;z) = Xz + \lrp{(-\eta-\sqrt{\theta})(\theta-1) X^3+X} z^3 \\
+ \lrp{X^9 + (\eta+\sqrt{\theta}(\theta-1)X^3} z^9 - X^9z^{27}.
\end{multline*}
Once again, the same formula holds $v$-adically. 


\subsection{\texorpdfstring{$\widetilde\sgn(F) = 1$}{tilde sgn of F is 1}, but \texorpdfstring{$\widetilde\sgn(F^{(1)})$}{tilde sgn of F twisted} transcendental}
\label{eg:sgnTranscendental}
\,\\
Back in Remark \ref{rmk:normalizationOfShtukaFunction}, we have promised an example of a function $F$ with $\widetilde\sgn(F) = 1$, but $\widetilde\sgn (F^{(1)})$ transcendental over $\bb F_q$. Here we provide such an example: suppose $X = \bb P^1$ with function field $\bb F_q(t)$, and set $\overline \infty$ to be corresponding to $t - c$ for some $c \in \overline{\bb F_q}$. Fix $\sgn$ such that $\sgn$ of the minimal polynomial of $\infty$ is 1. Then the function 
\[
F = \frac{c - c^q}{c - \theta} \frac{t-\theta}{t-c^q}
\]
has $\widetilde\sgn(F) = 1$. For $d_\infty > 2$, we have that 
\[
\widetilde{\sgn}(F^{(1)}) = \lrp{\frac{c-c^q}{c-\theta}}^q \frac{c - \theta^q}{c - c^{q^2}} = \frac{c^q - c^{q^2}}{c - c^{q^2}} \frac{c - \theta^q}{c^q - \theta^q},
\]
which is transcendental over $\bb F_q$. If $d_\infty = 2$, then from 
\[
\widetilde\sgn ((t - c)(t - c^q)) = 1,
\]
we have that 
\[
\widetilde\sgn (t - c) = \frac{1}{c - c^q}. 
\]
Thus 
\[
\widetilde{\sgn}(F^{(1)}) = \lrp{\frac{c-c^q}{c-\theta}}^q (c - \theta^q)(c - c^q) = (c-c^q)^{q+1} \frac{c - \theta^q}{c^q - \theta^q},
\]
which is also transcendental over $\bb F_q$. 


\subsection{\texorpdfstring{$X = \bb P^1$}{P1 (g = 0)}, any \texorpdfstring{$d_\infty$}{degree of infinity}}
\,\\
To end our list of examples, we will have a glimpse on how the $v$-adic convergence should work for a general curve, by studying Hayes modules coming from $\bb P^1$ other than the Carlitz module. 

We continue using the settings from the previous example \ref{eg:sgnTranscendental}. Let $X = \bb P^1$ with function field $\bb F_q(t)$. Set $\overline \infty$ to be the point corresponding to the place $t - c$ for some $c \in \overline{\bb F_q}$. This will only exclude the case with $\infty$ corresponding to the usual degree at $t$, which gives the well-understood Carlitz module. 

Once again, fix $\sgn$ such that $\sgn$ of the minimal polynomial of $\infty$ is 1. A shtuka function $f$ has divisor 
\[
\mathrm{div} (f) = V^{(1)} - V + (\Xi) - (\overline \infty^{(1)}),
\]
and is given by
\[
f = C \frac{t - \theta}{t - c^q},
\]
where $C \in \bb C_\infty$ is a constant. Since
\[
\widetilde\sgn (f f^{(1)} \cdots f^{(d_\infty - 1)}) = 1
\]
(see Remark \ref{rmk:normalizationOfShtukaFunction}), we can compute that
\[
C^{\frac{q^{d_\infty}-1}{q-1}} = \lrp{(c-\theta)(c-\theta^q) \cdots (c - \theta^{q^{d_\infty-1}})}^{-1}.
\]

We can also rewrite the shtuka functions as 
\[
f = U \frac{1}{c - \theta} \frac{t-\theta}{t-c^q},
\]
where 
\[
U^{\frac{q^{d_\infty}-1}{q-1}}
= \frac{(c - \theta)^{\frac{q^{d_\infty}-q}{q-1} } }{(c-\theta^q) \cdots (c - \theta^{q^{d_\infty-1}})}.
\]
The advantage of writing $f$ in this way is that as an element in $H$, the $\dfrac{q^{d_\infty}-1}{q-1}$-th power of $U$ only has nonzero valuation at places above $\infty$ (as in $K$, not $\bb K$). The notation $U$ comes from the fact that the $\dfrac{q^{d_\infty}-1}{q-1}$-th power of $U$ is a unit in $\mc O_{H}$. In particular, $U$ has valuation zero at the place above ``degree in $\theta$'', which is not true for $C$. This will help us to compute the $v$-adic convergence of $e_\rho$ and $\log_\rho$ for $v$ to be the place ``degree in $\theta$''.

Now fix such a $U$ (out of the choice of a $\dfrac{q^{d_\infty}-1}{q-1}$-th root of unity). The exponential series for the Hayes module corresponding to this choice is given by 
\begin{align*}
e_\rho(z)
= z &+ U^{-1} (c-\theta) \frac{(\theta - c)^q}{\theta^q-\theta} z^q \\
&+ U^{-(q+1)} (c-\theta)^{1+q} \frac{(\theta^q - c)^q (\theta - c)^{q^2}}{(\theta^{q^2} - \theta)(\theta^q-\theta)^q} z^{q^2} + \cdots \\
= z &- \sum_{n=1}^\infty U^{-\frac{q^n-1}{q-1}} (c-\theta)^{\frac{q^{n+1}-1}{q-1}} \\
&\phantom{=======} \cdot \frac{(\theta^{q^{n-1}} - c)^q (\theta^{q^{n-2}} - c)^{q^2} \cdots (\theta^q-c)^{q^{n-1}} }{(\theta^{q^n} - \theta) (\theta^{q^{n-1}} - \theta)^q \cdots (\theta^q-\theta)^{q^{n-1}}} z^{q^n}.
\end{align*}

The differential $\omega \in H^0(\overline X, \Omega^1(V - (\overline \infty) - (\overline \infty^{(-1)})))$ as in \cite[0.3.7]{Thakur1993_ShtukasAJ} and Remark \ref{rmk:NormalizationOfDifferential} is given by 
\[
\omega^{(1)} = -U \frac{dt}{(t-c^q)(t-c)},
\]
and the logarithm series is given by 
\begin{align*}
\log_\rho(z) 
= z &+ U^{-1} (c - \theta) \frac{(\theta - c)^q}{(\theta - \theta^q)} z^q \\
&+ U^{-(q+1)} (c-\theta)^{1+q} \frac{(\theta - c)^{q^2}(\theta - c^q)}{(\theta - \theta^q)(\theta - \theta^{q^2})} z^{q^2} + \cdots \\
= z &- \sum_{n = 1}^\infty U^{-\frac{q^n-1}{q-1}} (c - \theta)^{\frac{q^{n+1}-1}{q-1}} \\
&\phantom{=======} \cdot \frac{(\theta-c^q)(\theta-c^{q^2}) \cdots (\theta-c^{q^{n-1}})}{(\theta-\theta^q) (\theta - \theta^{q^2}) \cdots (\theta - \theta^{q^n})} z^{q^n}.
\end{align*}

From the calculation we have done in section \ref{sec:vForCoeff}, namely Corollary \ref{cor:sumOfvForISubk}, we can see that if $v$ is a place of $\bb P^1$ away from $\infty$ and the one corresponding to $\deg$ on $\theta$, then 
\begin{itemize}
\item $e_\rho(z)$ converges in $\bb C_v$ for all $w(z) > e_w \frac{1}{q^{\deg \mf p}-1}$, where $w, e_w$ are as in section \ref{sec:vForCoeff};
\item $\log_\rho(z)$ converges in $\bb C_v$ for all $v(z) > 0$. 
\end{itemize}
For $v$ corresponding to degree in $\theta$, we can directly compute the valuation. Let $w$ be a place in $H^+ := H(U)$ above $v$, and $e_w$ the ramification index of $w$ over $v$. Then
\[
w(e_n) = -e_w \frac{q^{n}-1}{q-1} > -e_w q^n \frac{1}{q-1},
\]
\[
w(l_n) = -e_w n,
\]
by looking at the degree in $\theta$ of the coefficients. This shows that the $v$-adic convergence behavior of $e_\rho(z)$ and $\log_\rho(z)$ is the same also when $v$ is the ``degree in $\theta$''. This gives evidence to the general $v$-adic convergence behavior of $e_\rho$ and $\log_\rho$ in the case when $d_\infty > 1$. 


%% file: Sec8.tex
Before we end this paper, we present two possible directions for generalizing results in this paper. 

\subsection{A conjecture for the general case}

From the previous $\bb P^1$ example and Theorems \ref{thm:vadicConvOfLog} and \ref{thm:vadicConvOfExp}, we formulate the following conjecture that predicts the general $v$-adic convergence for $e_\rho$ and $\log_\rho$. 

\begin{conj}
Let $\rho$ be a Hayes module on $X$, with no restriction on genus $X$ or $d_\infty$, and $H^+$ be the field of definition of $\rho$ (see \cite[\S7.3]{Goss1998_BasicStrcBook}, \cite[\S3.3]{Thakur2004}). Fix a place $v$ of $X$ different from $\infty$, normalized so that the value group is $\bb Z$. Let $K_V$ be the finite field extension of $H^+$ containing all zeros of the Drinfeld divisor $V$ corresponding to $\rho$, $w$ be the place (normalized so the value group is $\bb Z$) in $K_V$ over $v$ upon a fixed embedding $\overline K \to \overline K_v$, and $e_w$ be the ramification index of $w$ over $v$. Then:
\begin{enumerate}
\item the exponential series $e_\rho(z)$ converges in $\bb C_v$ when 
\[
w(z) > w(J_0) + e_w \cdot e_\theta \cdot \frac{1}{q^{\frac{\deg \mf p}{f_\theta}}-1},
\]
where $J_0$ is is an ideal depending on the Drinfeld divisor $V$ evaluated at $\Xi$, $e_\theta, f_\theta$ positive integers depending on ramification and intertia respectively;

\item the logarithm series $\log_\rho(z)$ converges in $\bb C_v$ when $v(z) > 0$. Moreover, the $v$-adic valuation of the coefficients of logarithm has order of magnitude $O(n)$. 
\end{enumerate}

\end{conj}

An immediate difficulty that we face when trying to prove this conjecture is that it is hard to explicitly write down 
\begin{enumerate} 
\item the shtuka function $f$ with a nice integral model, with all zeros being integral;
\item the differential $\omega$ in terms of the shtuka function $f$.
\end{enumerate}

In particular, the way we obtain the good presentation for the shtuka functions $f$ for elliptic curves and ramifying hyperelliptic curves is via long division, which requires that there is a degree 2 element in $\bb A$. Obviously this does not have to be true in general. For instance this will fail immediately when $d_\infty \geq 3$. In the general case, long division can still be done, but the result will no longer be in terms of integral functions, needless to say $\widetilde \sgn$ 1 functions. It also becomes hard to keep track of the relations of the generators in a model for $\bb A$, when there are more than 2 generators and/or more than 1 relation.


\subsection{Higher dimension}

It has been shown by Anderson and Thakur \cite[2.4.1]{Anderson_Thakur1990} that for the $n$-th tensor of Carlitz module $C^{\otimes n}$, the logarithm series $\mathrm{Log}_n(z)$ converges $v$-adically for all $z \in \bb C_v^n$ with $v(z) > 0$, and they used the convergence to calculate the values $\zeta_v(n)$. In \cite[Theorem 3.3.3]{ChangMishiba2017}, Chang and Mishiba showed a generalization for some uniformizable $t$-modules. 
As for elliptic curves with $\infty$ as the rational point at infinity, Green has given an expression of coefficients of the logarithm series for the $n$-tensor power of a Hayes module in his thesis \cite[Theorem 4.2.4]{GreenThesis2018}. A work in progress of the author is to prove a similar $v$-adic convergence result on the tensor product, which can be used to calculate $v$-adic zeta values. 


\subsection{Other \texorpdfstring{$L$}{L}-series}

It is also worth mentioning that Green and Papanikolas \cite{GreenPapanikolas2018} studied the shtuka functions for elliptic curves and come up with another formula for $L(1,\chi)$, as a special case for a formula for Pellarin $L$-series \cite{Pellarin2012}. It would be an exciting idea to see if we can come up with similar $v$-adic results for Pellarin $L$-series \cite{Pellarin2012}.